\input amstex
\UseAMSsymbols
\loadbold
\documentstyle{amsppt}
\magnification=\magstep1
\vsize=7.5in
\topmatter
\title{The prime number race and zeros
 of $L$-functions off the critical line, II}
\endtitle
\author Kevin Ford$^*$, Sergei Konyagin$^\dag$
\endauthor
\thanks
$^*$ First author supported by National Science Foundation grant
DMS-0070618.
\endthanks
\thanks
$^\dag$ Second author was supported by Grants 02-01-00248
and 00-15-96109 from the Russian Foundation for Basic Research.
\endthanks
\rightheadtext{Prime number race and zeros of $L$-functions, II}
\leftheadtext{K. Ford, S. Konyagin}
\date July 22, 2002 \enddate
\address\noindent 
K.F.: Department of Mathematics,
University of Illinois, Urbana, IL 61801
S.K.: Department of Mechanics and Mathematics, Moscow State University,
\newline
Moscow 111992, Russia.
\endaddress
\subjclass\nofrills{2000 Mathematics Subject Classification}
 Primary 11N13, 11M26 \endsubjclass
\abstract
 We continue our examination the effects of certain hypothetical
 configurations of zeros of Dirichlet $L$-functions lying off the critical
 line on the relative magnitude of the functions $\pi_{q,a}(x)$.
 Here $\pi_{q,a}(x)$ is
 the number of primes $\le x$ in the progression $a \mod q$.
 In particular, we look at situations where $\pi_{q,1}(x)$
 is simultaneously greater than,
 or simultaneously less than, each of $k$ functions
 $\pi_{q,a_i}(x)$  ($1\le i\le k$).  We also consider the total number of
 possible orderings of $r$ functions $\pi_{q,a_i}(x)$ ($1\le i\le r$).
\endabstract
\endtopmatter
\document

\predefine\barunder{\b}

\def\flr#1{\left\lfloor #1 \right\rfloor}

\def\sg{\sigma}
\def\g{\gamma}

\redefine\b{\beta}

\define\bc{\overline{\chi}}
\def\vp{\varphi}

\def\curly{\Cal}
\def\BB{\curly B}

\def\RR{\Bbb R}

\redefine\le{\leqslant}
\redefine\ge{\geqslant}

\def\a{\alpha}
\def\e{\varepsilon}
\define\({\left(}
\define\){\right)}
\define\pfrac#1#2{\( \frac{#1}{#2} \)}
\define\bfrac#1#2{\left[ \frac{#1}{#2} \right]}
\define\fq{F_q^*}
\define\sumstar{\sideset \and^* \to\sum}

\def\sg{\sigma}

\def\sgn{\text{sgn}}

\def\Arg{\operatorname{Arg}}
\def\FF{\curly F}
\def\GG{\curly G}

\def\II{\curly I}
\def\PP{\curly P}

\parindent 0.35in

\head 1.  Introduction \endhead

%
%
%
%

Denote by $\pi_{q,a}(x)$ the number of primes $p\le x$ with $p\equiv
a\pmod{q}$.
This paper is a continuation of our investigations from \cite{FK1}
on problems concerning the relative magnitude of
$\pi_{q,a}(x)$ for a fixed $q$ and varying $a$.
More about the background of the ``prime race'' problems may be found in
\cite{FK1} and \cite{FK2}.  As in \cite{FK1}
we are concerned with the consequences of hypothetical configurations of
zeros of Dirichlet $L$-functions lying off the critical line.
Roughly speaking, each zero of an $L$-function imparts an oscillation on
the functions $\pi_{q,a}(x)$, the zeros with largest real part giving the
largest oscillations.   In \cite{FK1} we were concerned with the
orderings of three functions $\pi_{q,a_i}(x)$ ($i=1,2,3$)
which occur for arbitrarily large $x$.
Let $\fq$ denote the multiplicative group of reduced residues modulo $q$.
 Our principal result, in
simple terms, was that for all $q\ge 5$ and distinct
$a_1,a_2,a_3\in \fq$, there are
finite configurations of hypothetical zeros which, if they really
existed, would imply that one of the orderings does not occur for large $x$.
Also, configurations can be constructed so the zeros all have imaginary
parts $\ge \tau$ for any given $\tau>0$.  The point of the exercise is this.
If one wishes to prove that all 6 orderings of the functions occur for
arbitrarily large $x$, one must prove in particular that our hypothetical
configurations are not possible.

In this paper we address two main types of problems.
First, if $D$ is a subset of $\fq\backslash \{1\}$, can it occur
for arbitrarily large $x$ that $\pi_{q,1}(x)$ is simultaneously smaller than,
or simultaneously large than,
each function $\pi_{q,a}(x)$ ($a\in D$) ?
Secondly, given
a subset $D$ of $\fq$, how many of the $|D|!$ possible orderings
of the functions $\pi_{q,a}(x)$ ($a\in D$) occur for arbitrarily large $x$ ?
In the language of Knapowski and Tur\'an, consider a game with players
$a_1, \ldots, a_k$, player $a_i$ having a score of $\pi_{q,a_i}(x)$ at time $x$.
Our questions can then be phrased as (i) Does player 1 lead infinitely
often or trail infinitely often? (ii) How many of the $|D|!$ orderings of the
players occur infinitely often?

Throughout, $q$ is a natural number, $q\ge 3$.
Below are some other definitions we will use.
$$
\split
C_q &= \text{ the set of non-principal Dirichlet characters modulo } q, \\
C_q(a,b) & = \{ \chi\in C_q: \chi(a) \ne \chi(b) \}, \\
\lambda(q) &= \text{ Carmichael's function: the largest order of an element
of }\fq, \\
\flr{x} &= \text{ the greatest integer which is } \le x, \\
\{ x\} &= x-\flr{x}, \text{ the fractional part of } x, \\
e(z) &= e^{2\pi i z}.
\endsplit
$$
Constants implied by the Landau $O-$ and
Vinogradov $\ll-$ symbols may depend on $q$, but not on any
other variable.

We begin with a lemma showing the relationship between functions
$\pi_{q,a}(x)$ and zeros of $L$-functions modulo $q$.

\proclaim{Lemma 1.1}
Let $q\ge 3$ and $a\in\fq$. Let $N_q(c)$ denote
the number of incongruent solutions of the congruence
$w^2\equiv c \pmod{q}$, and let $\pi(x)$ be the number of primes
$\le x$. Then for $x\ge 2$,
$$
\phi(q) \pi_{q,a}(x) = \pi(x) - 2 \Re \biggl( \sum_{\chi\in C_q}
\bc(a) \sumstar_{\Sb L(\rho,\chi)=0 \\ \Im \rho \ge 0 \\ \Re \rho > 0\endSb}
 f(\rho) \biggr)
 -N_q(a) \frac{x^{1/2}}{\log x} + O\(\frac{x^{1/2}}{\log^2 x}\),
$$
where
$$
f(\rho) := \frac{x^{\rho}}{\rho\log x} + \frac{1}{\rho}
\int_2^x \frac{t^{\rho}}{t\log^2 t}
\, dt= \frac{x^{\rho}}{\rho\log x} + O\pfrac{x^{\Re \rho}}{|\rho|^2\log^2 x},
$$
zeros are counted with mutiplicity,
and  $\sumstar$ indicates
that the summand is $\frac12 f(\rho)$ if $\Im \rho=0$.
\endproclaim

Lemma 1.1 is well-known, following from explicit formulas
(e.g. [Da], chapters 19, 20).  See also the proof of Lemma 1.1 of \cite{FK1}.

\proclaim{Corollary 1.2}
Let $\sg > \frac12$, $q\ge 3$ and $a,b\in\fq$.
Then as $x\to \infty$,
$$
\phi(q)\( \pi_{q,a}(x) - \pi_{q,b}(x) \) = -2\Re \biggl( \sum_{\chi\in C_q}
(\bc(a)-\bc(b)) \sumstar_{\Sb L(\rho,\chi)=0 \\ \Im \rho \ge 0 \\
\Re \rho \ge \sg \endSb}
 f(\rho) \biggr) + o\(\frac{x^\sg}{\log x}\).
$$
\endproclaim

Corollary 1.2 is a very old result, and follows from Lemma 1.1 and bounds
$$
\sum_{|\Im \rho| \ge x} \frac{x^\rho}{\rho} =o(x^{1/2}), \quad
N(T,\chi) \ll T\log T, \quad N(T,\sg,\chi)\ll_\sg T^{1-\delta(\sg)},
$$
where $\delta(\sg) > 0$ for $\sg>1/2$ and
$$
N(T,\chi)=|\{ \rho : |\Im \rho| \le T, \Re \rho > 0 \}|, \quad
N(T,\sg,\chi) =|\{ \rho : |\Im \rho| \le T, \Re \rho \ge \sg \}|.
$$
See for example a similar analysis for the approximation of $\pi(x)$ in
\cite{SP}.
The first two estimates above can be found in Davenport (\cite{Da}, Ch. 19, 20)
and an example of the third can be found in
Montgomery (e.g. [Mo], Theorem 12.1).  The upper bound on
$N(T,\sg,\chi)$ implies that
$$
\sum_{\chi\in C_q} \sum_{\Sb L(\rho,\chi)=0 \\ \Re \rho \ge \sg \\
\Im \rho > T\endSb} \frac{1}{|\rho|} \ll_\sg T^{-\delta(\sg)}. \tag{1.1}
$$

In applying Corollary 1.2, frequently we approximate $f(\rho)$ by
$x^\rho/(\rho \log x)$ with a total error of at most
$$
O \biggl(\frac{x^{\sg}}{\log^2 x} \sum_{\chi\in C_q(a,b)}
\sum_{\Sb L(\rho,\chi)=0 \\ \Re \rho \ge 1/2 \endSb}
\frac{1}{|\rho|^2} \biggr) = O \( \frac{x^{\sg}}{\log^2 x} \).
$$
Therefore we have the following.

\proclaim{Corollary 1.3}
Let $q\ge 3$, $a,b\in\fq$,
 $\sg > 1/2$ and suppose for $\chi\in C_q(a,b)$, the zeros of
$L(s,\chi)$ have real part $\le \sg$.  Then, as $u\to \infty$,
$$
\frac{u\phi(q)}{2e^{\sg u}} \( \pi_{q,a}(e^u) - \pi_{q,b}(e^u) \) =
\sum_{\chi\in C_q(a,b)} \;\;
 \sumstar_{\Sb L(\sg+it,\chi)=0 \\ t \ge 0 \endSb}
\frac{\nu(b)-\nu(a)}{\sqrt{t^2+\sg^2}}  + o(1),
$$
where $\nu(n) = \sin(t u - \Arg \chi(n) + \tan^{-1} (\sg/t))$.
Here we adopt the convention that $\tan^{-1}(\sg/t)=\pi/2$ if $t=0$.
\endproclaim
An inequality which is useful when $t$ is large is
$$
|\sin (v+\tan^{-1}(\sg/t)) - \sin(v) | \le \tan^{-1}(\sg/t) \le \sg/t.
\tag{1.2}
$$

Questions concerning the signs of the differences
$\pi_{q,a}(x) - \pi_{q,b}(x)$ therefore boil down to questions
about the trigonometric sums occurring in Lemma 1.1 and Corollaries 1.2,1.3.
As opposed to \cite{FK1}, a {\it barrier} in this paper refers
to the existence of a system of trigonometric sums of this type with
certain properties, and has nothing directly to do with prime counting
functions.  All of our results on the existence or
non-existence of particular types of barriers have consequences
for the distribution of
functions $\pi_{q,a}(x)$, but it is important to separate the
two.

Suppose for each $\chi\in C_q$, $B(\chi)$ is a sequence of complex numbers
with non-negative imaginary part (possibly empty, duplicates allowed), and
denote by $\BB$ the system of $B(\chi)$ for $\chi\in C_q$.  Let $n(\rho,\chi)$
be the number of occurrences of the number $\rho$ in $B(\chi)$.
If $\rho$ is real, we suppose that $n(\rho,\chi)=n(\rho,\bc)$.
The sets $B(\chi)$ will play the role of hypothetical zeros of the
$L$-function $L(s,\chi)$.  Define
$$
R^+(\BB) = \sup \{ \Re \rho: \rho \in \BB \}, \quad
R^-(\BB) = \inf \{ \Re \rho: \rho \in \BB \}.
$$
We shall suppose throughout that
$$
\frac12 < R^-(\BB) \le R^+(\BB) \le 1 \tag{1.3}
$$
and also, in accordance with (1.1), that
$$
\sum_{\chi\in C_q} \sum_{\rho\in B(\chi)} \frac{n(\rho,\chi)}{|\rho|} <
\infty. \tag{1.4}
$$
In accordance with Lemma 1.1, define
$$
P_{q,a}(x;\BB) = -\frac{2}{\phi(q)}
\Re \biggl( \sum_{\chi\in C_q} \bc(a)
\sumstar_{\rho\in B(\chi)} n(\rho,\chi) f(\rho) \biggr) + \frac{\pi(x)}
{\phi(q)} \tag{1.5}
$$
and
$$
D_{q,a,b}(x;\BB) = P_{q,a}(x;\BB) - P_{q,b}(x;\BB).
$$
where as before $\sumstar$ means the inner summand is $\frac{n(\rho,\chi)}2
 f(\rho)$ when
$\rho$ is real.  We say that two functions $F_1,F_2:[0,\infty)\to \RR$
are $\b$-similar if $|F_1(x)-F_2(x)| = o(x^\b/\log x)$ as $x\to
\infty$.  This is related to the conclusions in Lemma 1.1 and Corollary 1.2.
For indexed sets of functions $\FF=\{ F_i \}, \GG=\{ G_i \}$,
we say that $\FF$ and $\GG$ are $\b$-similar if $F_i$ and $G_i$
are $\b$-similar for each $i$.  With $q$ and $\BB$ fixed,
let $\PP_q$ be the list of functions $P_{q,a}(x;\BB)$.
For a system of functions $\FF$, also indexed by $a\in \fq$, suppose
$\II(\FF)$ is a statement concerning the
magnitudes of functions $F_{q,a}(x)$.
An example is
$$
\text{For sufficiently large $x$, at least one }
F_{q,a}(x) < F_{q,1}(x) \quad (a\in\fq\backslash\{1\}).
$$
For a system $\BB$, let $\b=R^-(\BB)$.
We say that $\BB$ is a {\it barrier for $\curly I$} if, for every
$\FF$ which is $\b$-similar to $\PP_q$, $\II(\FF)$ is false.

To relate this to the prime race problem,
let $\Pi_q$ be the list of functions $\pi_{q,a}(x)$,
indexed by $a\in\fq$.  Let $z_\BB$ denote the condition that for each
$\chi\in C_q$ and $\rho\in B(\chi)$,  $L(s,\chi)$ has a zero
of multiplicity $n(\rho,\chi)$ at $s=\rho$, and all other zeros of $L(s,\chi)$
in the upper half plane have real part less than $R^{-}(\BB)$.
By Lemma 1.1, if $z_\BB$ then $\Pi_q$ is $\b$-similar to $\PP_q$, thus
we have the following.

\proclaim{Lemma 1.4}  If $\BB$ is a barrier for $\II$ and
condition $z_\BB$ holds, then $\II(\Pi_q)$ is false.
\endproclaim

If each sequence $B(\chi)$ is finite, we call $\BB$ a {\it finite
barrier} for $\II$ and denote by $|\BB|$ the sum of the number of elements of
each sequence $B(\chi)$, counted according to multiplicity. We say that
$|\BB|$ is the size of the barrier $\BB$.
Of primary interest is to construct barriers for $\II$ where the
imaginary parts of the points in each $B(\chi)$ are all $\ge \tau$
for an arbitrarily large $\tau$.
It may occur that $|\BB|$ remains bounded as $\tau\to \infty$, in which
case we say that $\II$ possesses a {\it bounded barrier} (which is
actually a sequence of barriers).  Later we will demonstrate
the non-existence of bounded barriers for certain statements $\II$.
There is one more type of barrier which we will work with, the {\it extremal
barrier}, which will be defined in section 4.  Finally, we remark that
in general we can choose $R^-(\BB)$
and $R^+(\BB)$ arbitrarily as long as $1/2<R^-(\BB)$, $R^+(\BB)<1$.

\medskip

An important feature of the sums $D_{q,a,b}(x;\BB)$ is that the
``dominant parts'' are often almost periodic functions.  To be specific, let
$$
g(\rho)=g(\rho;a,b)=\sum_{\chi\in C_q} n(\rho,\chi) (\chi(a)-\chi(b)),
\quad \b(a,b)= \sup \{ \Re \rho : g(\rho)\ne 0 \}.
\tag{1.6}
$$
Also let
$$
z(\chi;a,b) = \{ \rho\in B(\chi) : g(\rho)\ne 0, \Re \rho = \b(a,b) \},
\quad z(a,b) = \bigcup_{\chi\in C_q} z(\chi;a,b).
\tag{1.7}
$$
In essence, the numbers in $z(a,b)$ are the ones which produce the
dominant terms in $D_{q,a,b}(x;\BB)$, provided $z(a,b)$ is non-empty.
Writing $\b=\b(a,b)$ for brevity, we have
$$
\split
D_{q,a,b}(x;\BB) &=  \frac{2x^{\b}}{\phi(q)\log x} M_{q,a,b}(x;\BB) +
 E_{q,a,b}(x;\BB), \\
M_{q,a,b}(x;\BB) &:= - \Re \biggl( \sum_{\chi\in C_q} (\bc(a)-\bc(b))
\sumstar_{\b+i\g\in z(\chi;a,b)} n(\b+i\g,\chi) \frac{x^{i\g}}{\b+i\g} \biggr)
\\
&= - \Re \biggl(\; \sumstar_{\b+i\g\in z(a,b)} \overline{g(\b+i\g)}
\frac{x^{i\g}}{\b+i\g} \biggr).
\endsplit\tag{1.8}
$$
Using Lemma 1.1 and (1.4), we have
$$
\split
|E_{q,a,b}(x;\BB)| &\ll \sum_{\chi\in C_q}
  \biggl[\,\,\; \sumstar_{\rho\in B(\chi)}
  \frac{x^{\Re \rho}}{|\rho|^2 \log^2 x} +
  \sumstar_{\Sb \rho\in B(\chi) \\ \rho\not\in z(\chi;a,b) \endSb}
  \frac{x^{\Re \rho}}{|\rho| \log x} \biggr] \\
&\ll \frac{x^{\b}}{\log^2 x} + o\( \frac{x^{\b}}{\log x}\) \\
&=o\( \frac{x^{\b}}{\log x} \)\quad(x\to\infty).
\endsplit\tag{1.9}
$$

A function $f$ is said to be {\it almost
periodic} with respect to a norm $\| \cdot \|$ if for and $\e>0$,
there is an $L>0$, so that any
real interval of length $L$ contains a number $\tau$ so that
$$
\| f(u+\tau) - f(u) \| \le \e.
$$
It follows from (1.4) and Theorems 8 and 12 of \S 1 of Chapter 1
in \cite{Be} that each sum $M_{q,a,b}(e^u;\BB)$
is a uniformly continuous almost periodic function in the sense of Bohr;
that is, almost periodic with respect to the supremum norm.
If one takes $B(\chi)$ to be the set of zeros $\rho$ of $L(s,\chi)$
with $\Re \rho=\b$ and $\Im \rho \ge 0$ (for $\chi\in C_q(a,b)$),
then $M_{q,a,b}(x;\BB)$ is precisely the double sum appearing in the
conclusion of Corollary 1.3 with $\sg=\b$.  Thus this double sum is also
a uniformly continuous almost periodic function in the sense of Bohr.
For a uniformly continuous almost periodic function $f$, define
$$
\|f\|_2=\lim_{U\to\infty}\left(\frac1U\int_0^U |f^2(u)|du\right)^{1/2}
$$
(the limit exists by Theorem 2 of \S 3 of Chapter 1 in \cite{Be}).
Next, if $f_1,\ldots,f_k$ are almost periodic with respect to a
norm $\| \cdot \|_A$, then the vector-valued function
$$
f(u) = (f_1(u),\ldots,f_k(u))
$$
is almost periodic with respect to the norm
$$
\| f \|_B := \max_{1\le j\le k} \| f_j \|_A.
$$

If, for some $\chi\in C_q$, $\chi(a)\ne \chi(b)$ and
 all non-trivial zeros of $L(s,\chi)$
have real part 1/2 (the Extended Riemann Hypothesis for $\chi$),
the inner sum in Corollary 1.2 (with $\sg=1/2$) is not uniformly convergent
(in fact, it has infinitely many jump discontinuities), but it is
still almost periodic in the sense of Stepanov (\cite{Be},
chapter 2).  That is, it is almost periodic with respect to the
norm
$$
\| g\|_{S^2} := \max_{x\in \Bbb R} \( \int_{x}^{x+1} |g(y)|^2\, dy
\)^{1/2}.
$$
The proof of this is implicit in \cite{K2}; another proof and generalization
can be found in \cite{KR}.  We note that if
a function is almost periodic in the Bohr sense, it is also almost
periodic in the Stepanov sense, since $\| g \|_{S^2} \le \| g\|_{\infty}$.
Any function $g$, almost period function in the Stepanov sense,
has the property that if $u$ is a continuity point of $g$, then
for every $\e>0$ there is an unbounded set of $v$ so that
$|g(v)-g(u)| \le \e$.

{\bf Remark 1.1.}
When each function in a set $\FF$ is almost periodic in the
Stepanov sense, to prove that some set of (strict) inequalities among a set
of functions $\FF$ occurs for an an bounded set of $u$, it suffices to prove
that the set of inequalities occur for a single $u$ which is a continuity
point of each function.  We can in fact make a stronger conclusion:
for some $L$ and $\delta>0$, on any interval of length $L$,
the measure of the set of $u$ for which the set of inequalities occur is
$\ge \delta$. As a consequence, setting $u=\log x$, we conclude that
the set of inequalities occurs on a set of $x$ of positive lower asymptotic
density.

{\bf Acknowledgement.}   Much of this paper was written while the authors
enjoyed the generous hospitality of the Mathematisches
 Forschungsinstitut Oberwolfach.

\vfil\eject
\head 2. Signs and comparison of trigonometric polynomials \endhead

First, we formulate some simple properties of trigonometric polynomials.
In particular, we prove that a real $n$-term trigonometric polynomial with 
a zero constant term must be nonnegative on a large set. By $\mu(E)$, where 
$E\subset\bold R$, we denote the Lebesgue measure of $E$. 
\proclaim{Lemma 2.1} Let $P$ be a real trigonometric  polynomial\newline 
$P(u)=\sum\sp n\sb {k=1}c\sb k\sin{(t\sb ku+\a_k)}$ 
$(c\sb k\in\bold R,t\sb k\in\bold R,t_k\neq0, t_k\neq t_l(k\neq l))$,
$E_+=\{u:\ P(u)\ge0\}$. Then\newline
1) $\int_0^U P(u)du=o(U)\quad(U\to\infty)$;\newline
2) $\|P\|_2:=\lim_{U\to\infty}
\(\frac1U\int_0^U(P(u))^2du\)^{1/2}
=\(\frac12\sum_{k=1}^nc_k^2\)^{1/2}$;\newline
3) $\|P\|_\infty\ge \|P\|_2\ge\sqrt{\frac1{2n}}\sum_{k=1}^n|c_k|$, where
$\|P\|_\infty:=\sup_u|P(u)|$;\newline
4) $\sup_u P(u)\ge\max_k |c_k|/2\ge\|P\|_{\infty}/2n$;\newline
5) $\mu(E_+\cap[0,U])/U\ge\frac1{4n}+o(1)\quad(U\to\infty)$.
\endproclaim
\demo{Proof} We have
$$\int_0^U P(u)du=\sum_{k=1}^n \frac{c_k}{t_k}(\cos\a_k-\cos(t_kU+\a_k)).
\tag2.1$$
The right-hand side of (2.1) is bounded for $u\in\bold R$, and 1) follows.

Further,
$$
\split
\int_0^U(P(u))^2du&=\sum_{k,l=1}^nc_kc_l\int_0^U  
\sin{(t_ku+\a_k)}\sin{(t_lu+\a_l)}du\\
&=\sum_{k=1}^n\int_0^U\(\frac {c_k^2}2
+\frac {c_k^2}2\cos((2t_ku+2\a_k))\)\\
&\quad +\sum_{k\neq l}^nc_kc_l\int_0^U  
\sin{(t_ku+\a_k)}\sin{(t_lu+\a_l)}du,
\endsplit
$$
and, by 1), 
$$\int_0^U(P(u))^2du=U\sum_{k=1}^n\frac{c_k^2}2+o(U)\quad(U\to\infty).$$

The first part in 3) follows from the inequality
$$\int_0^U(P(u))^2du\le\int_0^U\|P\|_\infty^2.$$ 
The second part follows from 2) and the Cauchy-Schwarz inequality.

To prove 4), we take $l$ so that $|c_l|=\max_k|c_k|$. Without loss of 
generality, $c_l\ge0$. Denote $a=\sup_u P(u)$. For $U>0$ we have
$$
\split
0&\ge\int_0^U(P(u)-a)(\sin{(t_lu+\a_l)
+1)}du \\ &=-aU-a\int_0^U \sin{(t_lu+\a_l)}du
+\int_0^U P(u)du+\int_0^U P(u)\sin{(t_lu+\a_l)}du \\ 
&=-aU +\int_0^U P(u)\sin{(t_lu+\a_l)}du+o(U)\quad(U\to\infty),
\endsplit\tag2.2
$$
by 1). Further, using again 1) and 2), we get
$$
\split
\int_0^U P(u)\sin{(t_lu+\a_l)}du
&=c_l\int_0^U (\sin{(t_lu+\a_l)})^2du\\
&\quad +\sum_{k\neq l}\int_0^U c_k\sin{(t_ku+\a_k)}\sin{(t_lu+\a_l)}du \\
&=\frac{c_l}2U+o(U)\quad(U\to\infty),
\endsplit\tag2.3
$$
and 4) follows from (2.2) and (2.3).

Denote $E_-=\{u:\ P(u)\le0\}$. By 1), 
$$\int_{E_+\cap[0,U]}P(u)du=-\int_{E_-\cap[0,U]}|P(u)|du
+o(U)\quad(U\to\infty).$$
Therefore,
$$\int_{E_+\cap[0,U]}P(u)du=\frac12\int_0^U|P(u)|du+o(U)
\quad(U\to\infty).\tag2.4$$
On the other hand, taking again $|c_l|=\max_k|c_k|$, we have, by (2.3),
$$\int_0^U|P(u)|du\ge\left|\int_0^UP(u)\sin{(t_lu+\a_l)}du\right| 
=\frac{|c_l|}2U+o(U)\quad(U\to\infty).$$
The equality (2.4) implies
$$\int_{E_+\cap[0,U]}P(u)du\ge\frac{|c_l|}4U+o(U)
=\frac{\max_k|c_k|}4U+o(U)\quad(U\to\infty).\tag2.5$$
Note that
$$\max_k|c_k|\ge\frac1n\sum_{k=1}^n|c_k|\ge\frac1n\|P\|_\infty.$$
Therefore,
$$\int_{E_+\cap[0,U]}P(u)du\le\|P\|_\infty\mu(E_+\cap[0,U])
\le n\max_k|c_k|\mu(E_+\cap[0,U]).\tag2.6$$
We exclude a trivial case when $P$ is not identically zero. 
Then combination of (2.5) and (2.6) proves 5) and thus  
completes the proof of Lemma.
\qed\enddemo

\proclaim{Theorem 2.2}\cite{N} Let $P$ be an exponential polynomial 
$P(u)=\sum\sp n\sb {k=1}c\sb ke\sp{it\sb ku}$ 
$(c\sb k\in\bold C,t\sb k\in\bold R)$, $U>0$,
$E\subset [0,U]$ of positive Lebesgue measure: 
Then
$$\max\sb {u\in [0,U]}\vert P(u)\vert
\leq\bigg\{\frac{CU}{\mu(E)}\bigg\}\sp {n-1}\sup\sb {u\in E}\vert P(u)\vert,$$
where $C$ is an absolute constant.
\endproclaim

\proclaim{Corollary 2.3} Let $P$ be a real trigonometric  polynomial\newline 
$P(u)=\sum\sp n\sb {k=1}c\sb k\sin{(t\sb ku+\a_k)}$ 
$(c\sb k\in\bold R,t\sb k\in\bold R,t_k\neq0, t_k\neq t_l(k\neq l))$,
$0<\g<1$, $S=\sum_{k=1}^n|c_k|$, $\e=\frac1{2\sqrt n}(C/\g)^{1-2n}$,
where $C$ is the constant from Theorem 2.2,
$E=\{u:\ |P(u)|<\e S\}$. Then for sufficiently large $U$
$$\mu(E\cap[0,U])/U<\g.$$
\endproclaim
\demo{Proof} Using Lemma 2.1, we get for sufficiently large $U$
$$\max_{u\in[0,U]}|P(u)|\ge\frac{S}{\sqrt{3n}}.$$
Suppose that 
$$\mu(E\cap[0,U])/U\ge\g.\tag2.7$$
Then, writing $P$ in an exponential form with $2n$ terms, we
get from Theorem 2.2,
$$\frac{S}{\sqrt{3n}}\leq(C/\g)\sp {2n-1}\sup\sb {u\in E}
\vert P(u)\vert\leq(C/\g)\sp {2n-1}\e S,$$
and, by the definition of $\e$, $S=0$, but in this case $E=\emptyset$. 
Thus, the supposition (2.7) cannot hold, and Corollary is proved. 
\qed\enddemo

\proclaim{Lemma 2.4} For any positive integer $n$ there exists such 
$\e_1=\e_1(n)>0$ that if $P,Q$ are real trigonometric polynomials,
$$P(u)=\sum\sp n\sb {k=1}a\sb k\cos{(t\sb ku+\a_k)},\quad 
Q(u)=\sum\sp n\sb {k=1}b\sb k\sin{(t\sb ku+\b_k)},$$ 
$$t_k\neq0,\ t_k\neq t_l(k\neq l),\quad 
|\a_k|\le\e_1,\ |\b_k|\le\e_1\quad(k=1,\dots,n),$$
then there exists a real number $u$ such that
$$P(u)\ge\e_1\sum_{k=1}^n|a_k|,\quad
Q(u)\ge\e_1\sum_{k=1}^n|b_k|.$$
\endproclaim

\demo{Proof} Take $\g=1/(10n)$. We will prove the lemma for
$$\e_1=\e/2,\tag2.8$$
where $\e$ is chosen in accordance with Corollary 2.3. Denote
$$\tilde P(u)=\sum\sp n\sb {k=1}a\sb k\cos{(t\sb ku)},\quad
S_1=\sum_{k=1}^n|a_k|.$$
Let $E=\{u:\ \tilde P(u)\ge0\}$, $E_1=\{u:\ |\tilde P(u)|<2\e_1S_1\}$.
Thus,
$$\forall u\in (E\setminus E_1)\quad \tilde P(u)\ge2\e_1S_1.\tag2.9$$

Take a sufficiently large $U$. By Lemma 2.1, we have 
$$\mu(E\cap[0,U])/U\ge\frac1{5n}.\tag2.10$$
Also, by Corollary 2.3 and (2.8),
$$\mu(E_1\cap[0,U])/U<\frac1{10n}.\tag2.11$$

Let 
$$\tilde Q(u)=\sum\sp n\sb {k=1}b\sb k\sin{(t\sb ku)},\quad
S_2=\sum_{k=1}^n|b_k|,\quad E_2=\{u:\ |\tilde Q(u)|<2\e_1S_2\}.$$
By Corollary 2.3 and (2.8),
$$\mu(E_2\cap[0,U])/U<\frac1{10n}.\tag2.12$$
The inequalities (2.10)---(2.12) show that the set
$E'=E\setminus E_1\setminus E_2$ is nonempty. Using 
evenness of $\tilde P$ and oddness of $\tilde Q$ we obtain
that for $u_1\in E'$ either $u=u_1$ or $u=-u_1$ satisfies the inequalities
$$\tilde P(u)\ge2\e_1S_1,\quad \tilde Q(u)\ge2\e_1S_2.$$
Taking into account, that, by the restrictions on $\a_k$ and $\b_k$,
we have
$$|P(u)-\tilde P(u)|\le\e_1S_1,\quad |Q(u)-\tilde Q(u)|\le\e_1S_2,$$
we get
$$P(u)\ge\e_1S_1,\quad Q(u)\ge\e_1S_2,$$
as required.
\qed\enddemo

The following lemma is closed to Lemma 1 from \cite{FFK}.
\proclaim{Lemma 2.5} Let $n$ be a positive integer, $0<\a<1$,
$$\e=\e(n,\a)=6(\a/6)^{2^{n-1}},$$ 
$s_1>\dots>s_n>0$. Then there exists a real number $u$ such that
$\e\le\{us_k\}\le\a$ for each $k\in\{1,\dots,n\}$. 
\endproclaim
\demo{Proof} We use induction on $n$. For $n=1$ we have $\e=\a$ and
the statement is trivial. Suppose that $n>1$ and the lemma holds for $n-1$.
We use the induction supposition for $\a'=\a^2/6$ instead of $\a$ and for
$\{s_2,\dots,s_n\}$. Observe that $\e=\e(n,\a)=\e(n-1,\a')$. There exists
a real number $u'$ such that $\e\le\{u's_k\}\le\a'$ for each 
$k\in\{2,\dots,n\}$. By Dirichlet's box principle, there exists a positive
integer $l$ satisfying $l\le3/\alpha$ and $\|lu's_1\|\le\alpha/3$. Take 
$u=lu'+\alpha/(2s_1)$. We have
$$\alpha/6\le\{us_1\}\le5\alpha/6$$
and for $k\in\{2,\dots,n\}$
$$\{us_k\}\le l\{u's_k\}+\alpha s_k/(2s_1)\le l\alpha'+\alpha/2\le\alpha,$$
$$\{us_k\}>\e,$$
as required.
\enddemo

\proclaim{Lemma 2.6} Let $n$ be a positive integer,  
$$\e_2=\e_2(n)=13^{-2^{n-1}},$$
$t_k$ be positive numbers, $|\b_k|\le\e_2$ for $k=1,\dots,n$.
Then there exists a real number $u$ such that  
$\sin(t_ku+\b_k)<-\e_2$ for each $k\in\{1,\dots,n\}$. 
\endproclaim
\demo{Proof} Take $\alpha=6/13$ and $s_k=t_k/(2\pi)$ for $k=1,\dots,n$.
By Lemma 2.5, there is $u'$ such that $\e_2\le\{u's_k\}\le\a$ for each 
$k\in\{1,\dots,n\}$. It is easy to check that $u=-2\pi u'$ 
satisfies Lemma 2.6.
\enddemo

\proclaim{Lemma 2.7} For any positive integer $n$ 
there exists such $\e_3=\e_3(n)>0$ that for any real $\gamma>0$
and real trigonometric polynomials
$$P(u)=\sum\sp n\sb {k=1}a\sb k\cos{(t\sb ku)},$$ 
$$Q(u)=\sum\sp n\sb {k=1}b\sb k\sin{(t\sb ku)},$$ 
$$R(u)=\sum\sp n\sb {k=1}c\sb k\sin{(t\sb ku)},$$ 
$$t_k\neq0, t_k\neq t_l(k\neq l),\quad 
b_k\ge|a_k|+c_k,\quad c_k\ge0,\quad(k=1,\dots,n),$$
$$\sum_{k=1}^n|a_k|>\gamma \sum_{k=1}^nb_k,$$ 
there exists a real number $u$ such that
$$Q(u)>\max(|P(u)|,R(u))+\e_3\gamma^2\sum_{k=1}^n b_k.$$
\endproclaim
The basic idea of the proof is the inequality 
$\|Q\|_2^2\ge\|P\|_2^2+\|R\|_2^2$ following from Lemma 2.1.
This inequality shows that there is a real $u$ such that 
$$Q^2(u)\ge P^2(u)+R^2(u)\ge\max(P^2(u),R^2(u)).\tag2.13$$
To strengthen (2.13), one can use the following possibilities:\newline
1) to estimate $\|Q\|_2^2-\|P\|_2^2-\|R\|_2^2$ from below;\newline
2) to estimate $\min(P^2(u),R^2(u))$ from below
and thus to strengthen the inequality
$P^2(u)+R^2(u)\ge\max(P^2(u),R^2(u))$;\newline
3) to show that $Q^2-P^2-R^2$ is not close to a constant
and thus has a big positive value at some point.\newline
It depends on the situation which of these arguments can work.
First we will prove a lemma using arguments 1) and 2).
\proclaim{Lemma 2.8} Under the suppositions of Lemma 2.7,
there exists $\e_4=\e_4(n)>0$ and a real number $u$ such that 
$$Q^2(u)>\max(P^2(u),R^2(u))+\max_k(\min((b_k^2-a_k^2-\gamma b_k^2)/2,
\e_4\gamma^2b_k^2)).$$
\endproclaim
\demo{Proof} Take any $k_0\in\{1,\dots,n\}$. If $c_{k_0}<\gamma b_{k_0}$,
then
$$\|Q\|_2^2-\|P\|_2^2-\|R\|_2^2=\frac12\sum_{k=1}^n(b_k^2-a_k^2-c_k^2)
\ge\frac12(b_{k_0}^2-a_{k_0}^2-c_{k_0}^2)>
\frac12(b_{k_0}^2-a_{k_0}^2-\gamma^2 b_{k_0}^2).$$
Therefore, there is $u$ such that
$$Q^2(u)-P^2(u)-R^2(u)>\frac12(b_{k_0}^2-a_{k_0}^2-\gamma^2 b_{k_0}^2).
\tag2.14$$
Now let us consider the case $c_{k_0}\ge\gamma b_{k_0}$, Let $\e$
be the number from Lemma 2.3, corresponding to $\gamma=1/3$,
$$E_1=\{u:\ |R(u)|<\e\gamma b_{k_0}\}.$$ 
By Lemma 2.3, for sufficiently large $U$
$$\mu(E_1\cap[0,U])/U<1/3.\tag2.15$$
Also, let
$$E_2=\{u:\ |P(u)|<\e\gamma b_{k_0}\}.$$ 
Taking into account the supposition of Lemma 2.7 for $\sum_{k=1}^n|a_k|$,
we get from Lemma 2.3
$$\mu(E_2\cap[0,U])/U<1/3.\tag2.16$$
Set $E_3=[0,U]\setminus E_1\setminus E_2$. By (2.15) and (2.16), we have
$$\mu(E_3)>U/3.\tag2.17$$
Also, from the definitions of $E_1$ and $E_2$ we find that for every 
$u\in E_3$
$$\min(|P(u)|,|R(u)|)\ge\e\gamma b_{k_0}.\tag2.18$$
Using Lemma 2.1, (2.17) and (2.18), we get
$$
\split
\int_0^U Q^2(u)du&\ge\int_0^U (P^2(u)+R^2(u))du+o(U)\\
&=\int_0^U \max(P^2(u),R^2(u))du
+\int_0^U \min(P^2(u),R^2(u))du+o(U)\\
&\ge\int_0^U \max(P^2(u),R^2(u))du
+\int_{E_3}(\e\gamma b_{k_0})^2du+o(U)\\
&\ge\int_0^U \max(P^2(u),R^2(u))du+(\e^2\gamma^2 b_{k_0}^2/3)U
+o(U)\quad(U\to\infty).
\endsplit
$$
Hence, there exists $u\in[0,U]$ such that
$$Q^2(u)-P^2(u)-R^2(u)>\e_4\gamma^2 b_{k_0}^2,\quad\e_4=\e^2/4.\tag2.19$$
So, for every $k_0$ one of the inequalities (2.14), (2.19) holds. This
proves Lemma 2.8.
\enddemo

\demo{Proof of Lemma 2.7} Without loss of generality, we can consider
$t_k>0$ for $k=1,\dots,n$ and $\gamma\le1/2$. Let
$\beta_1=\frac1n$, $\beta_j=\frac{\beta_{j-1}^2}{16}$ for $j=2,\dots,n+1$,
$$S=\sum_{k=1}^nb_k.$$
Choose the numbers $k_1,k_2,\dots$ so that
$$b_{k_1}\ge\beta_1S,\quad t_{k_j}>t_{k_{j-1}},\ b_{k_j}\ge\beta_jS
\ (j>1).$$
Note that $k_1$ can be always found because 
$$\max_kb_k\ge\frac1n\sum_kb_k=\beta_1S.$$
We terminate our construction when for some $l$ we cannot define a following 
number $k_{l+1}$, that is
$$\forall t_k>t_{k_l},\ b_k<\beta_{l+1}S.\tag2.20$$

If $|a_{k_l}|<b_{k_l}/2$, then, by Lemma 2.8,
there exists a real number $u$ such that 
$$Q^2(u)>\max(P^2(u),R^2(u))+\min((b_{k_l}^2-a_{k_l}^2
-\gamma b_{k_l}^2)/2,\e_4\gamma^2b_{k_l}^2).$$
Further,
$$b_{k_l}^2-a_{k_l}^2-\gamma b_{k_l}^2>\frac14b_{k_l}^2,$$
$$b_{k_l}\ge\b_nS.$$
Therefore,
$$Q^2(u)>\max(P^2(u),R^2(u))+\min(1/8,\e_4\gamma^2)\b_n^2S^2.\tag2.21$$

Now we have to consider the case 
$$|a_{k_l}|\ge b_{k_l}/2.\tag2.22$$
Define the even trigonometric polynomial $W(u)=Q^2(u)-P^2(u)-R^2(u)$
and estimate the coefficient $A$ of $\cos(Tu)$ in $W$, $T=2t_{k_l}$.
We have
$$\gather
A=-(a_{k_l}^2+b_{k_l}^2-c_{k_l}^2)/2-\sum\Sb t_{k}+t_{k'}=T,\\
t_{k_l}<t_k<2t_{k_l}\endSb (a_{k}a_{k'}+b_{k}b_{k'}-c_{k}c_{k'})\\
-\sum\Sb t_{k}-t_{k'}=T,\\t_k>2t_{k_l}\endSb 
(a_{k}a_{k'}+c_{k}c_{k'}-b_{k}b_{k'}).\tag2.23
\endgather$$
By (2.22),
$$(a_{k_l}^2+b_{k_l}^2-c_{k_l}^2)/2\ge
(a_{k_l}^2+b_{k_l}^2-(b_{k_l}-|a_{k_l}|)^2)/2=
|a_{k_l}|b_{k_l}\ge b_{k_l}^2/2\ge\beta_l^2S^2/2.\tag2.24$$
For $t_{k}>t_{k_l}$ and arbitrary $k'$ we have, by (2.20),
$$|a_{k}a_{k'}|+b_{k}b_{k'}+c_{k}c_{k'}
\le 2b_{k}b_{k'}\le2\beta_{l+1}b_{k'}S=\beta_l^2b_{k'}S/8.$$
Therefore,
$$\gather
\sum\Sb t_{k}+t_{k'}=T,\\
t_{k_l}<t_k<2t_{k_l}\endSb (a_{k}a_{k'}+b_{k}b_{k'}-c_{k}c_{k'})
+\sum\Sb t_{k}-t_{k'}=T,\\t_k>2t_{k_l}\endSb 
(a_{k}a_{k'}+c_{k}c_{k'}-b_{k}b_{k'})\\
\ge-2\sum_{k'}\beta_l^2b_{k'}S/8\ge-\beta_l^2S^2/4.\tag2.25
\endgather$$
Substituting (2.24) and (2.25) into (2.23), we get
$$A\le-\beta_l^2S^2/4\le-\b_n^2S^2/4.$$
By Lemma 2.1, taking into account that $W$ has a nonnegative constant term
we obtain
$$\sup_u W(u)\ge\b_n^2S^2/8.$$
Thus, in the case (2.22), there exists $u$ such that 
$$Q^2(u)-P^2(u)-R^2(u)>\b_n^2S^2/9.$$ 
In the opposite case we had the inequality (2.21).
So, for some $\e=\e(n)$ we always can find a real number $u_1$ such that
$$Q^2(u_1)>\max(P^2(u_1),R^2(u_1))+\e\gamma^2 S^2.$$ 
Let $x=|Q(u_1)|$, $y=\max(|P(u_1)|,|R(u_1)|)<x$. Using the inequality
$x-y>(x^2-y^2)/(2x)$ we get
$$\gather
|Q(u_1)|>\max(|P(u_1)|,|R(u_1)|)+\e\gamma^2 S^2/(2|Q(u_1)|)\\
\ge\max(|P(u_1)|,|R(u_1)|)+\e\gamma^2 S/2,
\endgather$$
and either $u=u_1$ or $u=-u_1$ satisfies the required inequalities
with $\e_4=\e/2$. Lemma 2.7 is proved.
\enddemo

\vfil\eject
\head 3.  Player 1 leading and trailing \endhead

%
%

\def\FAL{\text{FAL $x$}}
\def\FSL{\text{For sufficiently large $x$}}
For short, we abbreviate the phrase ``For arbitrarily large $x$'' by
``FAL $x$''.
In this section we address questions of whether or not
$$
\align
\FAL,\, \pi_{q,1}(x) & < \pi_{q,a}(x) \quad (\forall a\in D), \tag{3.1} \\
\FAL,\, \pi_{q,1}(x) & > \pi_{q,a}(x) \quad (\forall a\in D), \tag{3.2}
\endalign
$$
for various subsets $D$ of $\fq \backslash \{1\}$.
The residue $1 \mod q$ is special because it is the identity in $\fq$,
and this allows one to prove results about comparing $\pi_{q,1}(x)$ to
$\pi_{q,a}(x)$ which would be difficult otherwise.  For example,
in the cases $q=3,4,6$, $D=\{q-1\}$,
Littlewood \cite{Li} proved each of (3.1) and (3.2).
Knapowski and Tur\'an \cite{KT1} proved that under
the assumption that for each $\chi\in C_q$, $L(s,\chi)$ has
no zeros on the real segment $(0,1)$ (known as Haselgrove's condition
for $q$) that the difference $\pi_{q,1}(x)-\pi_{q,a}(x)$ changes sign
infinitely often.  Assuming the real parts of the nontrivial zeros
of $L(s,\chi)$ are all 1/2 for $\chi\in C_q$,
Kaczorowski \cite{K2} proved that
$$
\align
\FAL,\,
\pi_{q,1}(x) & < \pi_{q,a}(x) \quad
 (\forall a\in \fq \backslash \{1\}),\tag{3.3}\\
\FAL,\,
\pi_{q,1}(x) & > \pi_{q,a}(x) \quad
(\forall a\in \fq \backslash \{1\}). \tag{3.4}
\endalign
$$
In fact his proof gives
a little bit more: if $D\subset \fq, 1\not\in D$, and all nontrivial
zeros of $L(s,\chi)$ $(\chi \in \cup_{a\in D} C_q(a,1))$ have real part
$1/2$, then each of the inequalities (3.1) and (3.2) is true.

The statements pertaining to barriers which correspond to (3.1)--(3.4) are
$$
\align
\FAL,\, F_{q,1}(x) & < F_{q,a}(x) \quad (\forall a\in D), \tag{3.1'} \\
\FAL,\, F_{q,1}(x) & > F_{q,a}(x) \quad (\forall a\in D), \tag{3.2'} \\
\FAL,\,
F_{q,1}(x) & < F_{q,a}(x) \quad
(\forall a\in \fq \backslash \{1\}),\tag{3.3'}\\
\FAL,\,
F_{q,1}(x) & > F_{q,a}(x) \quad
 (\forall a\in \fq \backslash \{1\}). \tag{3.4'}
\endalign
$$
Among the results of this section we show the existence of bounded barriers for
(3.3') and (3.4') when $q\ge 7$, $q\not\in \{8,10,12,24\}$,
and show that no finite barriers exist for (3.3') when $q\in \{8,12,24\}$.
We also show that no bounded barriers exist for (3.3') and (3.4') when
$q\in \{5,10\}$.

For fixed $q$ define the quantities (analogs of (1.6), (1.7))
$$
\split
N(\rho,\chi) &= \text{ the multiplicity of the zero $\rho$ of } L(s,\chi), \\
G(\rho) &= G(\rho;a,b) = \sum_{\chi\in C_q} N(\rho,\chi) (\chi(a)-\chi(b)),\\
\sg(a,b) &= \sup \{ \Re \rho : G(\rho)\ne 0 \}, \\
Z(\chi;a,b) &= \{ \rho : L(\rho,\chi)=0, G(\rho)\ne 0,\Re \rho=\sg(a,b)\},\\
Z(a,b) &= \bigcup_{\chi\in C_q(a,b)} Z(\chi;a,b).
\endsplit
$$
The condition that $Z(a,b)$ is nonempty means that the supremum of the
real parts of the zeros $\rho$ of $L(s,\chi)$ with  $\chi\in C_q(a,b)$
and $G(\rho)\ne 0$ is attained.
In this case the sums over zeros in Corollary 1.3
are almost periodic functions
in the Stepanov sense.   In the case $b=1$, the condition
$G(\rho)\ne 0$ is equivalent to the statement that $L(\rho,\chi)=0$
for some $\chi$ with $\chi(a)\ne 1$ (in fact $\Re G(\rho)<0$ in this case).

%
%

\proclaim{Theorem 3.1}  Suppose $q\ge 3$, $D\subset F_q^*$ and $1\not\in D$.
Suppose $\BB$ is a system such that
for each $a\in D$ the set $z(a,1)$ is nonempty.
Then $\BB$ is a barrier for the statement $\II(\FF)$:
$$
\FSL,\, \exists a\in D : F_{q,1}(x) \ge \frac{F_{q,a}(x)+F_{q,a^{-1}}(x)}2.
$$
\endproclaim

\proclaim{Corollary 3.2} Suppose $q\ge 3$, $D\subset F_q^*$, $1\not\in D$,
and for each $a\in D$, $a^2\equiv 1\pmod{q}$.
If $\BB$ is a system such that $z(a,1)$ is nonempty for $a\in D$,
then (3.1') holds.
Consequently, there are no finite barriers for (3.3') when $q\in \{8,12,24\}$.
\endproclaim

\proclaim{Corollary 3.3}  Suppose $q\ge 3$, $D\subset F_q^*$ and $1\not\in D$.
If $Z(a,1)$ is non-empty for each $a\in D$, then
$$
\FAL,\, \pi_{q,1}(x) < \frac{\pi_{q,a}(x)+\pi_{q,a^{-1}}(x)}{2}
\quad (\forall a\in D).
$$
If in addition for each $a\in D$, $a^2\equiv 1\pmod{q}$, then (3.1) holds.
In particular, if $q\in\{8,12,24\}$ and $Z(a,1)$ is nonempty for
$a\in F_q^*\setminus\{1\}$, then (3.3) holds.
\endproclaim

\demo{Proof of Theorem 3.1}
We have $C_q(a,1)=C_q(a^{-1},1)$ and $z(a,1)=z(a^{-1},1)$ for $a\in D$.
For each $\chi\in C_q(a,1)$,
$(\bc(a) + \bc(a^{-1}))/2 - 1 = \Re \chi(a) - 1$ is a negative real number.
Let $\beta_a = \beta(a,1)$ for each $a\in D$ and put $\b=R^-(\BB)$.
Clearly $\b \le \min_{a\in D} \b_a$.  Let $\FF$ be $\b$-similar
to $\PP_q$.
By (1.8) and (1.9), for each $a\in D$ we have as $u\to\infty$
$$
\frac{u \phi(q)}{e^{\beta_a u}} \(\! F_{q,1}(e^u) -
\frac{F_{q,a}(e^u)+F_{q,a^{-1}}(e^u)}{2}\!\)
= - 2\!\! \sum_{\chi\in C_q(a,1)} (1-\Re \chi(a)) R_a(u;\chi) + o(1),
\tag{3.5}
$$
where
$$
R_a(u;\chi) = \frac{n(\beta_a,\chi)}{2\beta_a} +
\sum_{\Sb \gamma\in z(\chi;a,1) \\ \gamma>0 \endSb}
n(\b_a+i\g,\chi)
\frac{\sin(\g u + \tan^{-1} (\b_a/\g))}{\sqrt{\g^2+\b_a^2}}.\tag{3.6}
$$
Since each $z(a,1)$ is nonempty, it follows that
for each $a\in D$ one of the functions $R_a(u;\chi)$ is
not identically zero.
Each function $R_a(u;\chi)$ is almost periodic in the sense of Bohr.
To prove the theorem it suffices to show that there is
a $u$  for which each $R_a(u;\chi) > 0$ (among those functions which
are not identically zero).  Clearly $u=0$ is such a number.
\qed
\enddemo

\demo{Proof of Corollary 3.3}
Let $\sg_a=\sg(a,1)$ for $a\in D$, $A_1=\{ a : \sg_a > 1/2 \}$,
$A_2 = \{ a: \sg_a = 1/2 \}$, $\b=\min_{a\in A_1} \sg_a$.
For each $\chi\in C_q$, let $B(\chi)$ be the sequence
of all zeros of $L(s,\chi)$ with real part $\ge \b$, so $z_\BB$ holds.
If $A_2$ is empty, the Corollary follows from Lemma 1.4.
Otherwise, by Lemma 1.1 we have for each $a\in A_2$,
$$
\split
\frac{u \phi(q)}{e^{u/2}} \(\! \pi_{q,1}(e^u) -
\frac{\pi_{q,a}(e^u)+\pi_{q,a^{-1}}(e^u)}{2} \!\)
&= - 2\!\! \sum_{\chi\in C_q(a,1)} (1-\Re \chi(a)) R_a(u;\chi) \\
& \quad + (N_q(a)-N_q(1)) + o(1)\quad(u\to\infty).
\endsplit
$$
where
$$
R_a(u;\chi) = N(1/2,\chi) +
\sum_{\Sb L(1/2+it,\chi)=0 \\ t>0 \endSb}
N(1/2+it,\chi) \frac{\sin(t u + \tan^{-1} (1/2t))}{\sqrt{t^2+1/4}}.
$$
We always have $N_q(a) \le N_q(1)$.
As in the proof of Lemma 3.1, each $R_a(u;\chi)$ is positive in a
neighborhood of $u=0$ when $a\in A_1$.
When $a\in A_2$, $R_a(u;\chi)$ is continuous on $(0,\log 2)$
and $R_a(u;\chi)\to + \infty$ as $u\to 0^+$ (\cite{K1}; \cite{K2}, Lemma 2).
Therefore if $u$ is positive and sufficiently small,
it is a continuity point for all $R_a(u;\chi)$ and each $R_a(u;\chi)>0$.
\qed
\enddemo

The next results address inequalities (3.1'), (3.2') when $D$ is a cyclic
subgroup of $\fq$ or order 3.

%
%

\proclaim{Theorem 3.4}  Suppose $q\ge 3$ and $G=\{1,a,a^2\}\subset F_q^*$ is
a cyclic group of order 3. Suppose $\BB$ is a system such that
the set $z(a,1)$ is nonempty and
consists of numbers with imaginary part $\ge 2+\sqrt3$.
Then $\BB$ is a barrier for the statements:
$$
\FSL,\, F_{q,1}(x) \ge \min(F_{q,a}(x),F_{q,a^2}(x),
$$
$$
\FSL,\, F_{q,1}(x) \le \max(F_{q,a}(x),F_{q,a^2}(x),
$$
\endproclaim

\proclaim{Corollary 3.5}  Suppose $q\ge 3$ and $G=\{1,a,a^2\}\subset F_q^*$
is a cyclic group of order 3. If $Z(a,1)$ is non-empty and, in the case
$\sg(a,1) > 1/2$, $Z(a,1)$ consists of numbers with imaginary part
$\ge 2+\sqrt3$, then
$$
\FAL,\, \pi_{q,1}(x)<\min(\pi_{q,a}(x),\pi_{q,a^2}(x)),
$$
$$
\FAL,\, \pi_{q,1}(x)>\max(\pi_{q,a}(x),\pi_{q,a^2}(x)).
$$
\endproclaim

Corollary 3.5 can be deduced from Theorem 3.4 in the same way as we have
proved Corollary 3.3.

\demo{Proof of Theorem 3.4} Let $\beta = \beta(a,1)$ and put $\b_0=R^-(\BB)$.
Clearly $\b_0\le\b$.  Let $\FF$ be $\b_0$-similar
to $\PP_q$. For $j=1,2$ let $K_j=\{\chi\in C_q:\chi(a)=e(j/3)\}$.
By (1.8) and (1.9), we have
$$
\frac{\phi(q) u}{2e^{\b u}} \( F_{q,a^j}(e^u)-F_{q,1}(e^u) \) =
f(u) + (-1)^j g(u) + o(1),\quad(u\to\infty)
$$
where
$$
\split
f(v) &= \frac32 \sum_{\g\in z(a,1)} \frac{n(\g )}{\sqrt{\g^2+\b^2}}
\sin(\g v+\tan^{-1}\tfrac{\b}{\g})=f_1(v)+f_2(v),\\
f_1(v)&=\frac32 \sum_{\g\in z(a,1)}\frac{n(\g )}{\g^2+\b^2}\g\sin(\g v),\quad
f_2(v)=\frac32 \sum_{\g\in z(a,1)}\frac{n(\g )}{\g^2+\b^2}\b\cos(\g v),\\
n(\g ) &= \sum_{\chi\in K_1 \cup K_2} n(\b+i\g,\chi),\\
g(v) &= \frac{\sqrt3}2 \sum_{\g\in z(a,1)}\frac{m(\g )}{\sqrt{\g^2+\b^2}}
\cos(\g v+\tan^{-1} \tfrac{\b}{\g})=g_1(v)-g_2(v),\\
g_1(v)&= \frac{\sqrt3}2 \sum_{\g\in z(a,1)}\frac{m(\g )}{\g^2+\b^2}\g\cos \g v,
\quad g_2(v)= \frac{\sqrt3}2 \sum_{\g\in z(a,1)}\frac{m(\g )}{\g^2+\b^2}
\b\sin \g v,\\
m(\g ) &= \sum_{j=1}^2 (-1)^j \sum_{\chi\in K_j} n(\b+i\g,\chi).
\endsplit
$$
Since $\sum n(\g )/\sqrt{\g^2+\b^2}$ converges and $|m(\g )|\le n(\g )$,
the series in the definitions of the functions $f_1$, $f_2$, $g_1$, and $g_2$
are uniformly convergent, and thus these functions are Bohr almost periodic.
We need only find a single $v$ for which $f(v)-g(v)$ and $f(v)+g(v)$
are both positive, and a single $v$ for which $f(v)-g(v)$ and $f(v)+g(v)$
are both negative. Using the approximation of $f_1$, $f_2$, $g_1$,
and $g_2$ by trigonometric polynomials and Lemma 2.1,
$$
\split
\|\max( &|f_1|, |f_2|+|g_1|+|g_2|)-|f_2| -|g_1| -|g_2|\|_2\\
&\ge\|\max(|f_1|, |f_2|+|g_1|+|g_2|)\|_2-\|f_2\|_2 -\|g_1\|_2 -\|g_2\|_2\\
&\ge\|f_1\|_2-\|f_2\|_2 -\|g_1\|_2 -\|g_2\|_2 \\
&=\sqrt{\frac98\sum_{\g\in z(a,1)} \frac{\g^2n^2(\g )}{(\g^2+\b^2)^2}}
-\sqrt{\frac98\sum_{\g\in z(a,1)} \frac{\b^2n^2(\g )}{(\g^2+\b^2)^2}}
-\sqrt{\frac38\sum_{\g\in z(a,1)} \frac{\g^2m^2(\g )}{(\g^2+\b^2)^2}}\\
&\;\;
-\sqrt{\frac98\sum_{\g\in z(a,1)} \frac{\b^2m^2(\g )}{(\g^2+\b^2)^2}}
\ge\left(\sqrt{\frac98}-\sqrt{\frac38}\right)S_1-
\left(\sqrt{\frac98}+\sqrt{\frac38}\right)S_2,
\endsplit\tag{3.7}
$$
where
$$
S_1=\sqrt{\sum_{\g\in z(a,1)} \frac{\g^2n^2(\g )}{(\g^2+\b^2)^2}},\quad
S_2=\sqrt{\sum_{\g\in z(a,1)} \frac{\b^2n^2(\g )}{(\g^2+\b^2)^2}}.
$$
Further,
$$S_2/S_1\le\max_{\g\in z(a,1)}\b/\g<\max_{\g\in z(a,1)}1/\g\le1/(2+\sqrt3).$$
Substituting the last inequality into (3.7), we obtain
$$\gather
\|\max(|f_1|, |f_2|+|g_1|+|g_2|)-|f_2| -|g_1| -|g_2|\|_2\\
>\left(\sqrt{\frac98}-\sqrt{\frac38}-\frac1{2+\sqrt3}
\left(\sqrt{\frac98}+\sqrt{\frac38}\right)\right)S_1=0.
\endgather$$
Therefore, there exists $v_1$ such that
$$\max(|f_1(v_1)|, |f_2(v_1)|+|g_1(v_1)|+|g_2(v_1)|)
-|f_2(v_1)| -|g_1(v_1)| -|g_2(v_1)|>0,
$$
which is equivalent to\newline
$|f_1(v_1)|>|f_2(v_1)| +|g_1(v_1)| +|g_2(v_1)|$.
Observe that $f_1(-v_1)=-f_1(v_1)$, $|f_2(-v_1)|=|f_2(v_1)|$,
$|g_1(-v_1)|=|g_1(v_1)|$,$|g_2(-v_1)|=|g_2(v_1)|$. Thus, one of the numbers
$v\in\{v_1,-v_1\}$ satisfies the inequality
$$f_1(v)>|f_2(v)| +|g_1(v)| +|g_2(v)|,\tag3.8$$
and the other satisfies the inequality
$$-f_1(v)>|f_2(v)| +|g_1(v)| +|g_2(v)|.\tag3.9$$
The inequalities (3.8) and (3.9) imply $f(v)>|g(v)|$ and $f(v)<-|g(v)|$,
respectively, as required. This completes the proof of the theorem.
\qed
\enddemo

{\bf Remarks.}  R. Rumely \cite{R} has computed the small zeros of
$L$-functions modulo $q$ (with imaginary part
$\le 2600$) $3\le q\le 72$ and several larger $q$, and all such
zeros lie on the
critical line.  Thus for such $q$
the hypothesis in Corollary 3.5 about
the imaginary parts of the zeros in $Z_q(a,b)$ is satisfied.

%
%

Two following statements complement Theorem 3.1 and Corollaries 3.2 and 3.3
for the problem of winning.
\proclaim{Theorem 3.6}  For any $n$ there is an effectively computable number
$\tau$ such that if $q\ge3$, $D\subset F_q^*$, $1\not\in D$,
$\BB$ is a system such that for each $a\in D$ the set $z(a,1)$ is nonempty,
$\bigcup_{a\in D}z(a,1)$ consists of numbers with imaginary part $\ge\tau$
and contains at most $n$ elements, then $\BB$ is a barrier for the statement
$$
\FSL,\, \exists a\in D : F_{q,1}(x) \le \frac{F_{q,a}(x)+F_{q,a^{-1}}(x)}2.
$$
\endproclaim

\proclaim{Corollary 3.7} Suppose $q\ge 3$, $D\subset F_q^*$, $1\not\in D$,
and for each $a\in D$, $a^2\equiv 1\pmod{q}$. Then there are no bounded
barriers for (3.2'). Consequently, there are no bounded barriers for (3.4')
when $q\in \{8,12,24\}$.
\endproclaim

\proclaim{Corollary 3.8} For any $n$ there is an effectively computable
number $\tau$ such that if $q\ge3$, $D\subset F_q^*$,
$1\not\in D$, for each $a\in D$ we have $a=a^{-1}$,
for each $a\in D$ the set $Z(a,1)$ is nonempty, $\bigcup_{a\in D}Z(a,1)$
consists of numbers with imaginary part $\ge\tau$ and contains at most
$n$ elements, then (3.2) holds.
\endproclaim

\demo{Proof of Theorem 3.6}
We have $C_q(a,1)=C_q(a^{-1},1)$ and $z(a,1)=z(a^{-1},1)$ for $a\in D$.
For each $\chi\in C_q(a,1)$,
$(\bc(a) + \bc(a^{-1}))/2 - 1 = \Re \chi(a) - 1$ is a negative real number.
Let $\beta_a = \beta(a,1)$ for each $a\in D$ and put $\b=R^-(\BB)$.
Clearly $\b \le \min_{a\in D} \b_a$.  Let $\FF$ be $\b$-similar
to $\PP_q$. Take $\tau=1/\e_2$, where $\e_2=\e_2(n)$
was defined in Lemma 2.6. By (1.2) and Lemma 2.6,
there exists a real number $u$ such that
$\sin(\g u+\tan^{-1}(\b_a/\g))<-\e_2$ for each
$a\in D$ and $\g\in z(a,1)$. By periodicity of sines we can find
an arbitrary large $u$ satisfying these inequalities, and from (3.5)
and (3.6) (notice that under our suppositions
$n(\b_a,\chi)=0$) we deduce the assertion of the theorem.
\qed\enddemo


\proclaim{Theorem 3.9}  For any $n$ there is an effectively computable number
$\tau$ such that if $q\ge5$, $G\subset F_q^*$ is
a cyclic group of order 4, for each $a\in G\setminus\{1\}$ the set
$z(a,1)$ is nonempty, $\bigcup_{a\in G\setminus\{1\}} z(a,1)$
consists of numbers with imaginary part $\ge\tau$ and contains at most
$n$ elements, then $\BB$ is a barrier for the statements
$$
\FSL,\, \exists a\in G\setminus\{1\} : F_{q,1}(x) \ge F_{q,a}(x),
$$
$$
\FSL,\, \exists a\in G\setminus\{1\} : F_{q,1}(x) \le F_{q,a}(x).
$$
Consequently, there are no bounded barriers for (3.3') and (3.4')
when $q\in \{5,10\}$.
\endproclaim

\proclaim{Corollary  3.10} For any $n$ there is an effectively computable
number $\tau$ such that if $q\ge5$, $G\subset F_q^*$ is
a cyclic group of order 4, for each $a\in G\setminus\{1\}$ the set
$Z(a,1)$ is nonempty, $\bigcup_{a\in G\setminus\{1\}}Z(a,1)$
consists of numbers with imaginary part $\ge\tau$ and contains at most
$n$ elements, then (3.1) and (3.2) hold for $D=G\setminus\{1\}$.
\endproclaim

\demo{Proof of Theorem 3.9} Let $G=\{1,a_1,a_2,a_3\}$, $a_j=a_1^j$ for
$j=2,3$, $\beta_1 = \beta(a_1,1)$, $\b_2=\b(a_2,1)$ and $\b_0=R^-(\BB)$.
Clearly $\b_0\le\b_2\le\b_1$.  Let $\FF$ be $\b_0$-similar
to $\PP_q$. For $j=1,2,3$ let $K_j=\{\chi\in C_q:\chi(a_1)=e(j/4)\}$.
By (1.8) and (1.9), we have, as $u\to\infty$,
$$
\frac{\phi(q) u}{2e^{\b_1 u}} \( F_{q,a_1}(e^u)-F_{q,1}(e^u) \) =
f(u) + g(u) + o(1),\tag3.10
$$
$$
\frac{\phi(q) u}{2e^{\b_2 u}} \( F_{q,a_2}(e^u)-F_{q,1}(e^u) \) =
h(u) +  o(1),\tag3.11
$$
$$
\frac{\phi(q) u}{2e^{\b_1 u}} \( F_{q,a_3}(e^u)-F_{q,1}(e^u) \) =
f(u) - g(u) + o(1),\tag3.12
$$
where
$$
\split
f(v) &= \sum_{\g\in z(a_1,1)}
 \frac{k_1(\g)+2l(\g)}{\sqrt{\g^2+\b_1^2}} \sin(\g v+\tan^{-1}
\tfrac{\b_1}{\g}),\\
g(v) &= \sum_{\g\in z(a_1,1)}
\frac{m(\g)}{\sqrt{\g^2+\b_1^2}} \cos(\g v+\tan^{-1}\tfrac{\b_1}{\g}),\\
h(v) &= \sum_{\g\in z(a_2,1)}
\frac{2k_2(\g)}{\sqrt{\g^2+\b_2^2}} \sin(\g v+\tan^{-1}\tfrac{\b_2}{\g}),\\
k_j(\g) &= \sum_{\chi\in K_1 \cup K_3} n(\b_j+i\g,\chi) \quad (j=1,2),\\
l(\g) &= \sum_{\chi\in K_2} n(\b_1+i\g,\chi), \\
m(\g) &= \sum_{\chi\in K_1} n(\b_1+i\g,\chi)-\sum_{\chi\in K_3}
n(\b_1+i\g,\chi).
\endsplit
$$
Define $\e_2=\e_2(n)$ from Lemma 2.6.
We consider two cases.

Case I:  $\b_2 < \b_1$.  From (1.6) and (1.7), it follows that
$k_1(\g)=m(\g)=0$ for all $\g$.  By Lemma 2.6 and the almost periodicity
of $f(u)$ and $h(u)$, there are arbitrarily large $u$ so that
$$
f(u) < -2\e_2 \sum_{\g\in z(a_1,1)}
\frac{l(\g)}{\sqrt{\g^2+\b_1^2}}, \quad h(u) < -2\e_2
\sum_{\g\in z(a_2,1)}
\frac{k_2(\g)}{\sqrt{\g^2+\b_2^2}},
$$
whence (3.2') holds with $D=G\backslash\{1\}$.
Similarly, applying Lemma 2.6 to the functions $\sin(v-\tan^{-1} \b_j/\g)$,
there are arbitrarily large $u$ so that
$$
f(u) > 2\e_2 \sum_{\g\in z(a_1,1)}
\frac{l(\g)}{\sqrt{\g^2+\b_1^2}}, \quad h(u) > 2\e_2
\sum_{\g\in z(a_2,1)}
\frac{k_2(\g)}{\sqrt{\g^2+\b_2^2}},
$$
whence (3.1') holds.

Case II: $\b_2=\b_1$.  Write $\b=\b_1=\b_2$.
 Here we have $z(a_2,1)\subseteq z(a_1,1)=z(a_3,1)$ and $k_1(\g)=k_2(\g)$.
We again separate into two cases.

Case IIa: We have
$$
\sum_{\g\in z(a_1,1)}
\frac{|m(\g)|}{\sqrt{\g^2+\b^2}}\le
\frac{\e_2}2\sum_{\g\in z(a_1,1)}
\frac{k_1(\g)+2l(\g)}{\sqrt{\g^2+\b^2}}.\tag3.13
$$
By Lemma 2.6, there are arbitrarily large $u$ so that
$$
f(u) < -\e_2 \sum_{\g\in z(a_1,1)}
\frac{k_1(\g)+2l(\g)}{\sqrt{\g^2+\b_1^2}}, \quad h(u) < -2\e_2
\sum_{\g\in z(a_2,1)}
\frac{k_2(\g)}{\sqrt{\g^2+\b_2^2}}.
$$
Since $|m(\g)| \le k_1(\g)$, (3.13) implies that for such $u$,
$|g(u)| < \frac12 |f(u)|$.  Thus, by (3.10)--(3.12), (3.2') holds.
Similarly, applying Lemma 2.6 to the functions $\sin(v-\tan^{-1} \b/\g)$,
we see that (3.1') holds.

Case IIb: (3.13) does not hold.
By (3.10)--(3.12) and the almost periodicity of $f,g,h$,
the theorem will follow if we show that there are real $u$ and $v$ such that
$$
f(u)>\max(|g(u)|, f(u)-h(u)/2),\tag3.14
$$
$$
f(v)<\min(-|g(v)|, f(v)-h(v)/2).\tag3.15
$$
We approximate $f,g,f-h/2$ by the polynomials
$$
Q(u) = \sum_{\g\in z(a_1,1)}
 \frac{k_1(\g)+2l(\g)}{\sqrt{\g^2+\b^2}} \sin(\g u),
$$
$$
P(u) = \sum_{\g\in z(a_1,1)}
\frac{m(\g)}{\sqrt{\g^2+\b^2}} \cos(\g u),
$$
$$
R(u) = \sum_{\g\in z(a_1,1)}
 \frac{2l(\g)}{\sqrt{\g^2+\b^2}} \sin(\g u).
$$
Note that $|m(\g)|\le k_1(\g)$. Since (3.13) fails,
we can use Lemma 2.7 with $\gamma=\e_2/2$. Thus, there exists
a real number $u_0$ such that
$$
Q(u_0)>\max(|P(u_0)|,R(u_0))+\e S,\tag3.16
$$
where $\e=\e_3\gamma^2$, $S=\sum_\g\frac{k_1(\g)+2l(\g)}{\sqrt{\g^2+\b^2}}$.
The inequality (3.16) clearly implies
$$
Q(-u_0)<\min(-|P(-u_0)|,R(-u_0))-\e S.\tag3.17
$$
Taking into account (1.2), we get
$$
|f(u)-Q(u)|\le S/\tau,\quad|g(u)-P(u)|\le S/\tau,\quad
|f(u)-h(u)/2-R(u)|\le S/\tau.
$$
Therefore, we deduce (3.14) from (3.16) for $u=u_0$ and (3.15) from (3.17)
for $v=-u_0$ provided that $2S/\tau<\e S$, or $\tau>2/\e$. This completes
the proof.
\qed\enddemo

%
%

\proclaim{Theorem 3.11}  Let $q\ge 7$, $q\not\in \{8,10,12,24\}$.
There is a set $D\in\fq\backslash \{1\}$ with $|D|=3$ so that
for any $\tau>0$ there is a system $\BB$ with
$|\BB| \le 34$ which is a barrier for both inequalities (3.1')
and (3.2'),
and each sequence $B(\chi)$ consists of numbers with imaginary part $>\tau$;
\endproclaim

\demo{Proof}
The argument depends on the group structure of $\fq$.
Denote by $Z_k$ the cyclic group of order $k$.
Every $\fq$, $q\ge 7$,  $q\not\in \{8,10,12,24\}$, either contains a
cyclic group of even order $n\ge 6$ or contains a subgroup isomorphic
to $Z_4 \times Z_2$.  Our constructions depend on properties of
the functions
$$
\split
Q(v) &= 2 \sin v + \frac12 \sin(6v), \\
P(v) &= 2 \cos v - \frac12 \cos(6v), \\
R(v) &= \sum_{k=2}^7 \frac{p_k}{k} \sin(kv), \quad p_2=1, p_3=2, p_4=3,
p_5=4, p_6=3, p_7=2.
\endsplit
$$
The critical properties are
$$
\split
|P(v)| &> \sqrt{3} Q(v) \quad (0\le v\le 0.759, 2.7 \le v\le 2\pi), \\
R(v) &< 0 \quad (0.758 \le v < \pi).
\endsplit\tag{3.18}
$$
\medskip

We first consider the case when
$\fq$ has an element of even order $n$ with $n\ge 6$.
Without loss of generality, if $n$ is a power of 2, assume $n=8$.
Let $\chi$ be a character of order $n$,
and let $a$ be an element of $\fq$ of order $n$ such that
$\chi(a)=e(-1/n)$.
Fix $\beta>\frac12$ and large $\g>0$, let
$n(\b+ik\g,\chi^j)=m_{j,k}$
($1\le j\le n-1$, $1\le k\le K$).
Suppose $n(\rho,\chi)=0$ for all other pairs $(\rho,\chi)$.
Let $\FF$ be $\b$-similar to $\PP_q$. By (1.8) and (1.9), as $u\to\infty$,
$$
\frac{\phi(q) u\g}{2e^{\b u}} \( F_{q,a^r}(e^u)-F_{q,1}(e^u) \) =
G_0(u\g) - G_r(u\g) + O\( \frac{1}{\g} + \frac{1}{u} \) + o(1),
$$
where
$$
\split
G_r(v) &= \sum_{j,k} \frac{m_{j,k}}{k} \sin\(kv+\frac{2\pi jr}{n} \) \\
&=\sum_{j,k} \frac{m_{j,k}}{k} \left[ \cos \( \frac{2\pi jr}{n} \) \sin (kv)
+ \sin \( \frac{2\pi jr}{n} \) \cos(kv) \right].
\endsplit
$$
We take $D=\{a^s,a^{n-s},a^{n/2}\}$ for some $s\ne n/2$.  The theorem will
follow if we show that for every $v\in [0,2\pi)$, there
is a $r\in \{s,n-s,n/2\}$ so that $G_r(v)>G_0(v)$,
since $G_r(v)>G_0(v)$ implies $G_{n-r}(-v) > G_0(-v)$.

First, if $n=2^d h$, where $h$ is odd and $h\ge 3$, we take $m_{2,6}=3$,
$m_{2h-2,1}=2$, $m_{h,k}=p_k$ for $2\le k\le 7$ and $m_{j,k}=0$ for other
$j,k$, so $|\BB|=20$.  We obtain
$$
\split
G_0(v)-G_{2^d}(v) &= (1-\cos (4\pi/h)) Q(v) + \sin (4\pi/h) P(v), \\
G_0(v)-G_{n-2^d}(v) &= (1-\cos (4\pi/h)) Q(v) - \sin (4\pi/h) P(v), \\
G_0(v)-G_{n/2}(v) &= 2R(v).
\endsplit
$$
The theorem follows in this case from (3.22) and the fact that
$$
\left| \frac{1-\cos (4\pi/h)}{\sin (4\pi/h)} \right| \le \sqrt{3}.
$$
Next, suppose $n=8$ and take $m_{2,1}=4$, $m_{3,k}=m_{5,k}=p_k$ for
$2\le k\le 7$ and $m_{j,k}=0$ for other $j,k$, so $|\BB|=34$.  Then
$$
\split
G_0(v) - G_3(v) &= 4(\sin v - \cos v) + (2-\sqrt{2}) R(v), \\
G_0(v) - G_5(v) &= 4(\sin v + \cos v) + (2-\sqrt{2}) R(v), \\
G_0(v) - G_4(v) &= 4R(v).
\endsplit
$$
When $0\le v\le 0.758$ or $\pi \le v \le 2\pi$, one of the first two functions
is negative.

The last case is when
$\fq$ has a subgroup $G$ isomorphic to $Z_4 \times Z_2$.  Let
$\{a,b\}$ generate $G$, $a$ having order 4 and $b$ having order 2.
Let $\chi_1$ have order 4, $\chi_2$ have order 2 so that
$$
\chi_1(a)=-i, \chi_1(b)=1, \quad \chi_2(a)=1, \chi_2(b)=-1.
$$

Fix $\beta>\frac12$ and large $\g>0$, and let, for some $L$,
$n(\b+il\g,\chi_1^j\chi_2^k)=m_{j,k,l}$
for $0\le j\le 3, 0\le k\le 1, (j,k)\ne (0,0)$, $1\le l \le L$.
Suppose $n(\rho,\chi)=0$ for all other pairs $(\rho,\chi)$.
Let $\FF$ be $\b$-similar to $\PP_q$. By (1.8) and (1.9), as $u\to\infty$,
$$
\frac{\phi(q) u \g}{2e^{\b u}} \( F_{q,a^r b^s}(e^u)-F_{q,1}(e^u) \) =
G_{0,0}(u\g) - G_{r,s}(u\g) +  O\( \frac{1}{\g} + \frac{1}{u} \) + o(1),
$$
where
$$
G_{r,s}(v) = \sum_{j,k,l} \frac{m_{j,k,l}}{l} \sin\(lv + \frac{\pi}{2}rj + \pi
sk \).
$$
Note that $G_{0,0}(v) < G_{r,s}(v)$ implies $G_{0,0}(-v) > G_{4-r,2-s}(-v)$.
We take $D=\{a,a^3,b\}$.
Thus, if for all $v\in [0,2\pi)$,  $G_{0,0}(v) < G_{r,s}(v)$
for some pair $(r,s)\in \{ (1,0), (3,0), (0,1) \}$,
then $\BB$ is a barrier for both (3.1') and (3.2').
We take $m_{1,0,1}=1$ ($L(s,\chi_1)$ has a simple zero at $s=\b+it$),
and $m_{0,1,l} = p_l$ for $2\le l\le 7$,
Take $m_{j,k,l}=0$ for other $(j,k,l)$, so $|\BB|=16$.  Then
$$
\split
G_{0,0}(v) - G_{1,0}(v) &= \sin v - \cos v, \\
G_{0,0}(v) - G_{3,0}(v) &= \sin v + \cos v, \\
G_{0,0}(v) - G_{0,1}(v) &= 2R(v).
\endsplit
$$
When $0\le v < \pi/4$ or $3\pi/4 < v \le 2\pi$, we have $|\cos v| > \sin v$,
whence either $G_{0,0}(v) - G_{1,0}(v) < 0$ or
$G_{0,0}(v) - G_{3,0}(v) < 0$.  For the remaining $v$, $R(v)<0$ by (3.18).
\qed\enddemo

\proclaim{Corollary 3.12}  Let $q=5$ or $q\ge 7$.  Each inequality (3.3'),
(3.4') possesses a bounded barrier if and only if
$q\not\in \{5,8,10,12,24\}$.
\endproclaim

\vfil\eject
\head 4.  Extremal Barriers \endhead

%
%

By an {\it ordering} of the functions $\pi_{q,a_i}(x)$ $(1\le i\le r)$
we mean a chain of inequalities
$$
\pi_{q,a_{i(1)}}(x)\ge\pi_{q,a_{i(2)}}(x)\ge\dots\ge\pi_{q,a_{i(r)}}(x),
$$
where $\{i(1),\dots,i(r)\}$ is a permutation of $\{1,\dots,r\}$. Thus, we
admit non-strict inequalities in orderings, and in the case of coincidence
of some functions $\pi_{q,a_i}(x)$ several orderings occur for $x$.
Let $S_q(D)$ be the number of orderings of the functions
$\pi_{q,a}(x)$ ($a\in D$) which occur for arbitrarily large $x$.
Likewise, for a system $\BB$ and set of functions $\FF$,
define $s(D)=s(D;\FF)$ to be the number of orderings of functions
$F_{q,a}(x;\BB)$ ($a\in D$) which occur for arbitrarily large $x$.

If $\pi_{q,a}(x) > \pi_{q,b}(x)$ and $\pi_{q,a}(y) < \pi_{q,b}(y)$,
then $\pi_{q,a}(w)=\pi_{q,b}(w)$ at some point $w$ between $x$ and $y$.
This property of these
functions is crucial to results about $S_q(D)$.
If a set of functions $\FF$ has the property that for $f_i, f_j\in \FF$,
$f_i(x)<f_j(x)$ and
$f_i(y)>f_j(y)$ implies $f_i(w)=f_j(w)$ for some $w$ between $x$ and $y$,
we say that $\FF$ is {\it good}.

Let $D \subseteq \fq$ and $\b=R^-(\BB)$.
We say that $\BB$ is a {\it KT-system} (Knapowski-Tur\'an system)
for $D$, if for each set of functions $\FF$ which is
$\b$-similar to $\PP_q$ and every distinct $a,b\in D$,
$$
\FAL, \, F_{q,a}(x) > F_{q,b}(x).
$$
If $\BB$ is a KT-system for $D$ and $z_\BB$ holds, then
each difference $\pi_{q,a}(x)-\pi_{q,b}(x)$, $a,b\in D$, changes
sign infinitely often.  For several moduli $q$ this is known
unconditionally for all differences $\pi_{q,a}(x)-\pi_{q,b}(x)$,
$a,b\in F_q^*$, $a\neq b$ (see \cite{FK2}).
A KT-system $D$ has the property that
for distinct
$a,b\in D$ there is some $\rho$ with $g(\rho;a,b)\ne 0$,
for otherwise $D_{q,a,b}(x)$ is identically zero and one could take
$F_{q,c}(x)=P_{q,c}(x;\BB)$ for each $c\in D$.
In the opposite direction we have the following.

\proclaim{Proposition 4.1}
Let $D\subseteq \fq$.  Every system $\BB$ which lacks
real elements and for which $z(a,b)$ is nonempty for $a,b\in D$
is a KT-system for $D$.
\endproclaim

\demo{Proof} Take distinct $a,b\in D$, let $\b=R^-(\BB)$
and suppose $\FF$ is $\b$-similar to $\PP_q$.
By (1.5), (1.6), (1.8) and (1.9), as  $u\to\infty$
$$
\frac{u\phi(q)}{2e^{\b u}}\( F_{q,a}(e^u) - F_{q,b}(e^u) \)
= -h(u)+o(1),\tag4.1
$$
where $\b=\b(a,b)$ and
$$
h(u)=\sum_{\b+i\g\in z(a,b)}\Re \left(\overline{g(\b+i\g)}
\frac{e^{i\g u}}{\b+i\g}\right).\tag4.2
$$
By (1.4), the partial sums of (4.2) uniformly converge to $h$. By Lemma 2.1,
$$
\align
\sup_u  h(u) &\ge\sup_\g \frac{|g(\b+i\g)|}{4|\b+i\g|}, \\
\sup_u -h(u) &\ge\sup_\g \frac{|g(\b+i\g)|}{4|\b+i\g|}.
\endalign
$$
Taking into account that $h$ is almost periodic function in the Bohr sense,
we get from (4.1)
$$\liminf_{x\to\infty}\frac{\log x}{x^{\b}}
\( F_{q,a}(x) - F_{q,b}(x) \)<0,$$
$$\limsup_{x\to\infty}\frac{\log x}{x^{\b}}
\( F_{q,a}(x) - F_{q,b}(x) \)>0,$$
and the proposition is proved.
\qed\enddemo

\proclaim{Theorem 4.2}
If $\BB$ is a KT-system for $D=\{a_1,a_2,\dots,a_r\}$,
then for every good $\FF$ which is $\b$-similar to $\PP_q$ ($\b=R^-(\BB)$),
at least $r(r-1)/2+1$ orderings of the functions $F_{q,a_i}(x)$
occur for arbitrarily large $x$.
Consequently, under the condition $z_\BB$, $S_q(D) \ge r(r-1)/2+1$.
\endproclaim

\demo{Proof} Fix a good $\FF$ which is $\b$-similar to $\PP_q$.
Let us construct a graph $G$.  For each permutation
$P=(i(1),\dots,i(r))$ of the set $\{1,\dots,r\}$, let $N(P)$ be the set of
real $x\ge 1$ with
$$
F_{q,a_{i(1)}}(x)\ge F_{q,a_{i(2)}}(x)\ge\dots\ge F_{q,a_{i(r)}}(x).
$$
For each unbounded set $N(P)$, associate a vertex $v(P)$ of $G$.
Put an edge from $v(P_1)$ to $v(P_2)$ whenever (i) $P_2$ is obtained from
$P_1$ by transposing two neighbor elements $k,l$, and (ii) $N(P_1)\cap
N(P_2)$ is unbounded (note $x\in N(P_1)\cap N(P_2)$ implies
$F_{q,k}(x)=F_{q,l}(x)$).  Label such an edge by $\{k,l\}$.

Also, as $\BB$ is KT-system and $\FF$ is good,
for any numbers $i$ and $j$, $1\le i<j\le r$, there is an edge labeled
by $(i,j)$. We claim that the graph $G$ contains a subgraph $G'$ such that
each component of $G'$ is a tree and the
labelings of the edges in $G'$ contain again all
possible pairs $(i,j)$. Indeed if $G$ contains a cycle $H$ take two vertices
$g_1$ and $g_2$ from $H$. Then there are numbers $i$ and $j$ occurring in $g_1$
and $g_2$ in opposite orders. This means that in both arcs of the cycle $H$
connecting $g_1$ and $g_2$ there is an edge labeled by $(i,j)$. Delete one
of them.  We can repeat this procedure as long as the remaining graph contains
at least one cycle. In the end we get a required subgraph $G'$. The number of
edges of $G'$ is at least the number of distinct labels, thus it is at least
$r(r-1)/2$. Therefore, the number of vertices of $G'$ is
$\ge r(r-1)/2+1$.
\qed\enddemo

A system $\BB$ is called an {\it extremal barrier} for $D$
if it is a $KT$-system for $D$ and a barrier for the statement
$$
s(D) \ge \frac{r(r-1)}{2}+2.
$$
By Lemma 1.4, if $\BB$ is an
 extremal barrier and $z_\BB$ holds, at most $r(r-1)/2+1$ orderings of
the functions $\pi_{q,a}(x)$ ($a\in D$) occur for large $x$.  An
interesting problem is
to describe for each $q$ the sets $D$ possessing finite
extremal barriers. We are very far from a complete solution to this problem;
in particular, there is no $q$, $\vp(q)>2$, for which we know whether the
whole system $F_q^*$ has a finite extremal barrier.
In this section we present some results on existence and nonexistence of
extremal barriers. In particular we shall see that for large moduli $q$
there is a finite extremal barrier for some set $D$ with $|D|=r(q)\to\infty$
as $q\to\infty$.

\proclaim{Theorem 4.3}
For every cyclic group $G\subset F_q^*$ of order $r\ge6$
and for every set $D\subset G$ such that $1\not\in D$ and
$a^{-1}\not\in D$ if $a^{-1}\neq a\in D$, there is a bounded
extremal barrier for $D$.
\endproclaim
\proclaim{Remark}  The size of $\BB$ in our construction depends only on $r$,
and it can be effectively computed.
\endproclaim
To prove Theorem 4.3, we take a generator $a_1$ of the group $G$
and a character $\chi_1$ so that $\chi_1(a_1)=e(-1/r)$. For $j=1,\dots,r-1$
denote $a_j=a_1^j$, $\chi_j=\chi_1^j$. Take $\b_1\in (1/2,1)$,
large $\gamma$ and large positive integer $K$ depending on $r$. The idea
is to put $n(\b_1+ki\gamma,\chi_j)=N_{k,j}$
($k=1,\dots,K$, $j=1,\dots,r-1$), where $N_{k,j}$ are appropriate nonnegative
integers.
For $k=1,\dots,K$, $v=0,\dots,r-1$ define the functions
$$
G_{k,v}(u)=\sum_{j=1}^{r-1}\frac{N_{k,j}}k\sin(ku+2\pi jv/r).
$$
If $\FF$ is $\b$-similar to $\PP_q$, we have for $1\le v,w<r$
$$
\gather
F_{q,a_v}(x)-F_{q,a_w}(x)=\frac{2x^{\b_1}}{\gamma\log x}\times\\
\left(\sum_{k=1}^K G_{k,v}(\gamma\log x)-G_{k,w}(\gamma\log x)
+O\left(\frac1\gamma\right)\right) + o(1)\quad(x\to\infty).\tag4.3
\endgather
$$
To choose multiplicities $N_{k,j}$ we need the following lemma.

\proclaim{Lemma 4.4} Let $c_v$, $d_v$ ($v=0,\dots,r-1$)
be real numbers such that $c_v=c_{r-v}$ ($v=1,\dots,r-1$),
$d_0=0$, $d_v=-d_{r-v}$ ($v=1,\dots,r-1$). Then there exist
real numbers $\nu_j$ ($j=0,\dots,r-1$) such that
$$
\sum_{j=0}^{r-1}\nu_j\sin(u+2\pi jv/r)
=c_v\sin u+d_v\cos u\quad (v=0,\dots,r-1).\tag4.4
$$
\endproclaim
\demo{Proof} The system (4.4) is equivalent to the system of two systems of
linear equations
$$\sum_{j=0}^{[r/2]}\mu_j\cos(2\pi jv/r)=c_v
\quad(v=0,\dots,[r/2]),\tag4.5
$$
$$\sum_{j=1}^{[(r-1)/2]}\lambda_j\sin(2\pi jv/r)=d_v
\quad(v=1,\dots,[(r-1)/2]).\tag4.6
$$
where $\mu_0=\nu_0$, $\mu_j=\nu_j+\nu_{r-j}$, $\lambda_j=\nu_j-\nu_{r-j}$
($1\le j\le(r-1)/2$), $\nu_{r/2}=\mu_{r/2}$ if $r$ is even.
To prove solubility of the system (4.5) it suffices to check that
the system has no nontrivial solutions for $c_v=0$
$(v=1,\dots,[r/2])$. Assume the contrary. Consider the trigonometric
polynomial
$$
T(u)=\sum_{j=0}^{[r/2]}\mu_j\cos(2\pi ju).
$$
If not all $\mu_j$ are zero, the polynomial $T$ has at most $2[r/2]$
zeros on $[0,2\pi)$ counting with multiplicity. On the other hand,
by (4.5) with our supposition $c_v=0$, the points $2\pi v/r$
($v=0,\dots,r-1$) are zeros of $T$, and, moreover, $0$ is a double zero.
Hence, the total number of the zeros of $T$ on $[0,2\pi)$ counting with
multiplicity is at least $r+1>2[r/2]$. This contradiction shows that
$T\equiv0$. So, the system (4.5) has a unique solution for any $c_v$.
In the same way we can prove the solubility of the system (4.6).

Now we have the existence of numbers $\mu_j$ and $\lambda_j$
satisfying (4.5) and (4.6). To complete the proof of Lemma 4.4,
it remains to set $\nu_{r/2}=\mu_{r/2}$ for even $r$, $\nu_0=\mu_0$,
$\nu_j=(\mu_j+\lambda_j)/2$ for $1\le j<r/2$,
$\nu_j=(\mu_{r-j}-\lambda_{r-j})/2$ for $r/2<j<r$.
\qed\enddemo

Here we shall apply Lemma 4.4 for the case $c_v=0$ ($v=0,\dots,r-1$).
We have stated it for arbitrary $c_v$ taking into account other applications.

Let $V=\{v:a_v\in D\}$. Let us take
a system of continuous even $2\pi$-periodic functions $f_v$, $v\in V$,
and let us require the following properties to hold:\newline
1) If $r/2\in V$ then $f_{r/2}\equiv0$;\newline
2) $\int_0^\pi f_v(u)du=0$ for all $v\in V$;\newline
3) for every distinct $v\in V$ and $w\in V$, $v\neq w$, there is the
unique point $u=u_{v,w}\in[0,\pi]$ at which $f_v(u)=f_w(u)$, and,
moreover, for distinct (nonordered) pairs $(v,w)$ the points $u_{v,w}$
are distinct.\newline
Clearly, a system $\Omega=\{f_v\}$ exists; for example we can take several
functions in a general position from the set of piecewise linear functions
with zero average and one corner on $(0,\pi)$, with slope $0$ to the
right of $0$ and slope $1$ to the left of $\pi$.

Observe that the ordering of the functions $\{f_v(u)\}$, $u\in[0,\pi]$,
changes only at points $u_{v,w}$. On the other hand, if points
$u_1\in[0,\pi]$ and $u_2\in[0,\pi]$ are separated by some point $u_{v,w}$,
then $(f_v(u_1)-f_w(u_1))(f_v(u_2)-f_w(u_2))<0$. Thus, the number of
orderings of the functions $\{f_v(u)\}$, $u\in[0,\pi]$, is $|V|(|V|-1)/2+1$.
Since the functions $f_v$ are even and $2\pi-periodic$, this is the number
of orderings on the whole real line.

Take $U$ which is a multiple of $2\pi$ and arrange all the points
$\pm u_{v,w}+2\pi k\in(U,\infty)$ in increasing order $U<u_1<u_2<\dots$.
Denote $u'_0=U$, $u'_j=(u_j+u_{j+1})/2$ for $j\ge1$. We say that a system
$\tilde \Omega$ of continuous real functions $\tilde f_v(u)$,
$u\in[U,\infty]$, is of
$\Omega$-type if for any $j\ge0$ and for any $u\in[u'_j,u'_{j+1}]$ the
ordering of the functions $\{\tilde f_v(u)\}$ coincides with the ordering
$\{f_v(u'_j)\}$ or with the ordering $\{f_v(u'_{j+1})\}$ (recall that in the
case of some equalities $\tilde f_v(u)=\tilde f_w(u)$ we assign to the point
$u$ several orderings). The system $\Omega$ is an example of a system
of $\Omega$-type. We will repeatedly use the following simple fact.

\proclaim{Proposition 4.5} If a system $\tilde\Omega$ of $2\pi$-periodic
functions is of $\Omega$-type, then every system sufficiently
close to $\tilde\Omega$
in the uniform metric is of $\Omega$-type. Moreover, every system whose
pairwise differences are close to corresponding differences for
$\tilde\Omega$ in the uniform metric is of $\Omega$-type.
\endproclaim

Eventually, we shall show that the system $\{F_{q,a_v}(\gamma\log x)\}$,
$v\in V$, is of $\Omega$-type, which proves Theorem 4.3.

First, take any admissible system $\Omega=\{f_v:v\in V\}$.
We approximate the functions $f_v$ by even trigonometric polynomials
$T_v$ with zero average in the uniform norm. In the case $r/2\in V$
we take $T_{r/2}\equiv0$. Let
$$
T_v(u)=\sum_{k=1}^Kb_{k,v}\cos(ku).
$$
By Proposition 4.5, for sufficiently large $K=K(r)$ we can make the
approximation so good that the system $\{T_v\}$ is of $\Omega$-type.

By the conditions on $D$, $r-v\not\in V$ if $v\in V$ and $v\neq r/2$.
Let $b_{k,r-v}=b_{k,v}$ for $v\in V$, and set $b_{k,v}=0$ for
$v\not\in V$ and $r-v\not\in V$. By Lemma 4.4, there exist
real numbers $\nu_{k,j}$ ($k=1,\dots,K$, $j=0,\dots,r-1$) such that
$$
\sum_{j=0}^{r-1}\nu_{k,j}\sin(ku+2\pi jv/r)
=b_{k,v}\cos ku\quad (k=1,\dots,K; v=0,\dots,r-1).
$$
Therefore,
$$
T_v(u)=\sum_{k=1}^K\sum_{j=0}^{r-1}\nu_{k,j}\sin(ku+2\pi jv/r)
\quad (v\in V).
$$

Take a positive integer $N$ and define trigonometric polynomials
$$\tilde T_v(u)=\sum_{k=1}^K\sum_{j=0}^{r-1}\frac{\tilde N_{k,j}}{kN}
\sin(ku+2\pi jv/r)\quad (v\in V),$$
where
$$\tilde N_{k,j}=k[N\nu_{k,j}]\quad(k=1,\dots,K; j=1,\dots,r-1).$$
By Proposition 4.5, the system $\{\tilde T_v\}$ is of $\Omega$-type
provided that $N$ is large enough.

Finally, take $\tilde N=\min_{k,j}\tilde N_{k,j}$,
$N_{k,j}=\tilde N_{k,j}-\tilde N\ge0$. Since $1\not\in D$, we have
$0\not\in V$ and hence,
$$
\sum_{j=0}^{r-1}\sin(ku+2\pi jv/r)=0\quad(k=1,\dots,K; v\in V)
$$
and
$$
N\tilde T_v(u)=\sum_{k=1}^K G_{k,v}(u)\quad(v\in V).
$$
The equality (4.3) can be rewritten for $v,w\in V$ as
$$
\frac{\gamma\log x}{2Nx^{\b_1}}(F_{q,a_v}(x)-F_{q,a_w}(x))=
\tilde T_v(\gamma\log x)-\tilde T_w(\gamma\log x)
+O\left(\frac1\gamma\right)+o(1)\quad(x\to\infty).
$$
By Proposition 4.5, the system
$\{\frac{\gamma\log x}{2Nx^{\b_1}}F_{q,a_v}(x)\}$, $v\in V$,
is of $\Omega$-type on $[U,\infty)$ if $U$ and $\gamma$
are large enough. So is the system $\{F_{q,a_v}(x)\}$, $v\in V$,
as required.\qed
\medskip

It is not difficult to see that $\lambda(q)\to\infty$
as $q\to\infty$. A lower estimate
$$\lambda(q)>(\log q)^{c\log\log\log (q+20)}$$
with some $c>0$ was established in \cite{EPS}. Thus, we have the following.

\proclaim{Corollary 4.6} For sufficiently large $q$ there is
a finite extremal barrier for some set $D$ with
$|D|=r(q)\ge\lambda(q)/2\to\infty$ as $q\to\infty$.
\endproclaim

It is naturally to ask if there are bounded extremal barriers
for $D=G$. We show
that it is not so in the case $|G|=3$. However, we cannot
prove that for $|G|=3$ there are no finite extremal barriers.

\proclaim{Theorem 4.7} For any $n$ there is an effectively computable
number $\tau$ such that the following holds.
Let $q\in\bold N$, $a\in F_q^*$, $a^3=1$,
$G=\{1,a,a^2\}$. If $\BB$ is a system such that
$Z_q(a,1)$ and $Z_q(a,a^2)$ are nonempty and $Z_q(a,1)\cup Z_q(a,a^2)$
consists of numbers with imaginary part $\ge\tau$ and contains at most
$n$ elements, then $\BB$ is a barrier for the statement
$$
\FF\quad\text{is good and}\quad s(\{1,a,a^2\};\FF) \le 4.
$$
Consequently, under the condition $z_\BB$, $S_q(G) \ge 5$.
\endproclaim

To prove Theorem 4.7, we first
estimate the number of orderings if each of three players leads and
trails for arbitrarily large $x$.

\proclaim{Lemma 4.8} Let $D=\{a_1,a_2,a_3\}\subset F_q^*$ and $\BB$ be such
a system that for any function system $\FF$ which is $\b$-similar to $\PP_q$
($\b=R^-(\BB)$), and for any
$a'\in D$ there are arbitrary large $x$ and $y$ such that
$$F_{q,a'}(x)>\max(F_{q,a''}(x):\ a''\in D\setminus\{a'\}),$$
$$F_{q,a'}(y)<\min(F_{q,a''}(y):\ a''\in D\setminus\{a'\}).$$
Then for any good function system $\FF$ which is $\b$-similar to $\PP_q$,
at least 5 orderings of the functions $\{F_{q,a'}(x):a'\in D\}$
occur for arbitrary large $x$.
\endproclaim

\demo{Proof} By Theorem 4.2, since $\BB$ is a KT-system for $D$,
at least 4 orderings
occur for arbitrary large $x$.  Assume only 4 orderings occur.  Since $\FF$
is good, there are
arbitrary large $x$ such that $F_{q,a_1}(x)\ge F_{q,a_2}(x)=F_{q,a_3}(x)$ or
$F_{q,a_1}(x)\le F_{q,a_2}(x)=F_{q,a_3}(x)$.
Thus, in both cases for large $x$ there are at least 3 orderings where
$a_1$ leads or trails, and, therefore, at most one ordering where
$a_1$ is in the second position. The same holds for $a_2$ and $a_3$.
Hence, the number of orderings where some player is in the second position,
which is clearly the number of all orderings, is at most 3, but that
is impossible. Lemma 4.8 is proved.
\enddemo

\demo{Proof of Theorem 4.7} Let $\b_0=R^-(\BB)$ and a system $\FF$ be
good and $\b_0$-similar to $\PP_q$. Re-denote by $\tau'$ the number $\tau$
from Theorem 3.6. Take
$$
\tau=\max(\tau',1/\e_1(n)),
$$
where $\e_1(n)$ is the number from Lemma 2.4, and suppose that the
conclusion of the theorem does not hold. Then, by Lemma 4.8, one of the
players $1,a,a^2$ does not lead $\FAL$ or does not trail
$\FAL$.
We see from Theorem 3.4 that this is not player $1$.
Without loss of generality, assume that $a^2$ does not lead $\FAL$.
Thus, the orderings
$$
F_{q,a^2}(x)> F_{q,a}(x)\ge F_{q,1}(x),\quad
F_{q,a^2}(x)> F_{q,1}(x)\ge F_{q,a}(x)
$$
do not occur for large $x$.

We use notation and relationships form the proof of Proposition 4.1
 with $b=a^2$.
Note that for any $\chi\in C_q$ we have $\chi(a^2)=\overline{\chi(a)}$. Thus,
$g(\rho)=g(\rho;a,a^2)$ is a purely imaginary number, and
$$
\Re \left(\overline{g(\b+i\g)}\frac{e^{i\g u}}{\b+i\g}\right)
=\frac{\overline{g(\b+i\g)}}{i\sqrt{\g^2+\b^2}}
\cos(\g u+\tan^{-1}(\b/\g)).\tag4.7
$$
Now let us follow again the proof of Proposition 4.1 to approximate
$ P_{q,1}(e^u)- ( P_{q,a}(e^u)+ P_{q,a^2}(e^u))/2$. Let
$$
g_1(\rho)=\sum_{\chi}n(\rho,\chi)(1-(\chi(a)+\chi(a^2))/2),
$$
$$
\b_1=\max\{\Re\rho:\ g_1(\rho)\neq0\},
$$
$$
\Cal R_1=\{\rho:\ \Re\rho=\b_1,\ g_1(\rho)\neq0\},
$$
$$
h_1(u)=\sum_{\b_1+i\g\in\Cal R_1}\Re \left(\overline{g_1(\b_1+i\g)}
\frac{e^{i\g u}}{\b_1+i\g}\right).
$$
The formula (4.1) written for $(1,a)$ and $(1,a^2)$ gives
$$
\frac{u\phi(q)}{2e^{\b_1 u}}\(  P_{q,1}(e^u)-
( P_{q,a}(e^u)+ P_{q,a^2}(e^u))/2\)=-h_1(u)+o(1)\quad(u\to\infty).\tag4.8
$$
Now, $g_1(\rho)$ is always a real number, and therefore
$$
\Re \left(\overline{g_1(\b_1+i\g)}\frac{e^{i\g u}}{\b_1+i\g}\right)
=\frac{g(\b_1+i\g)}{\sqrt{\g^2+\b_1^2}}
\sin(\g u+\tan^{-1}(\b_1/\g)).\tag4.9
$$
Note, that in the definitions of $h$ and $h_1$ the sum is taken over
$\g\in Z_q(a,a^2)$ and, respectively, over $\g\in Z_q(a,1)$. By the choice of
$\tau$ and (1.2), any $\g\in Z_q(a,1)\cup Z_q(a,a^2)$ satisfies the
inequalities
$\tan^{-1}(\b/\g))<\e_1$, $\tan^{-1}(\b_1/\g))<\e_1$.
By the suppositions of the theorem, $h$ and $h_1$
are nonzero polynomials with at most $n$ distinct frequences $\g$ in total.
Therefore, we can apply Lemma 2.4 to $h$ and $h_1$. Hence, there exist
$\delta>0$ and $u\in\bold R$ such that
$$
h(u)>\delta,\quad h_1(u)>\delta.\tag4.10
$$
As the functions $h$ and $h_1$ are almost periodic in the Bohr
sense, we can find an arbitrary large $u$ satisfying (4.10). Then, by (4.1)
and (4.8), and the
$\b_0$-similarity of the system $\Cal F$ to $\PP_q$, taking into
account that $\b_0\le\min(\b,\b_1)$, we get
$$ F_{q,a^2}(x)>  F_{q,a}(x),\quad ( F_{q,a}(x)+  F_{q,a^2}(x))/2
>  F_{q,1}(x).$$
This contradicts our assumption that $a^2$ does not lead $\FAL$
and completes the proof of Theorem 4.7.
\qed\enddemo
\bigskip

\vfil\eject
\head 5.  The number of possible orderings \endhead

%
%
%
%

\proclaim{Theorem 5.1}
Fix $q$ and an arbitrarily large $\tau$.
There is a system $\BB$ satisfying
\item{(i)} $|\BB|$ bounded in terms of $q$;
\item{(ii)} $\rho\in \BB$ implies $\Im \rho > \tau$;
\item{(iii)} For every $r\ge 2$ distinct elements $a_1,\ldots,a_r$ of $\fq$, 
$\BB$ is a barrier for the property $s(\{a_1,\ldots,a_r\}) > r(r-1)$;
\item{(iv)} If $z_\BB$ holds, then for every $r\ge 2$ distinct elements
$a_1,\ldots,a_r$ of $\fq$, \newline $S_q(\{a_1,\ldots,a_r\}) \le r(r-1)$.
\endproclaim

\demo{Proof} 
Suppose $\fq$ is generated by $g_1,\ldots,g_m$,
which have orders $n_1,\ldots,n_m$, where $n_1 n_2 \cdots n_m = \phi(q)$.
Define $\chi_j$ by
$$
\bc_j(g_j) = e(1/n_j), \qquad \bc_j(g_h) = 1 \quad (h\ne j).
$$
Let $\g$ be large depending on $q$, and
$$
\tfrac12 < \b_m < \b_{m-1} < \cdots < \b_1 < 1.
$$
For $1\le j \le m, 1\le k\le 2$,
let $n(\b_j+ik\g,\chi_j^k) = c_{j,k}$.  Also, for each $j$ there is
at least one $k$ so that $n(\b_j+ik\g,\chi_j^k) \ge 1$.
In what follows, implied constants depend on $q, \g$ and the numbers 
$c_{j,k}$.  For each $a\in \fq$ write
$$
a \equiv g_1^{\a_1(a)} \cdots g_m^{\a_m(a)} \pmod{q}, \quad 0\le \a_j(a) \le
n_j-1.
$$
Let $\FF$ be $\b_m$-similar to $\PP_q$.  By (1.8) and (1.9),
$$
\split
\Delta_{a,b}(u) &:= -\frac{u\phi(q)}{2} \bigl[ F_{q,a}(e^u)-F_{q,b}(e^u) 
\bigr] \\
&= \sum_{j=1}^m e^{\b_j u} \( f_j(u,\a_j(a))-f_j(u,\a_j(b)) \)
+ o(e^{\b_m u}),\quad(u\to\infty)
\endsplit\tag{5.1}
$$
where
$$
f_j(u,\a)=\Re \sum_{k=1}^2 \frac{c_{j,k} e(k\a/n_j)}{\b_j+ik\g} \left[ 
e^{ik\g u} + u e^{-\b_j u} \int_2^{e^u}
\frac{v^{\b_j+ik\g}}{v\log^2 v}\, dv \right].
\tag{5.2}
$$
Let $J(a,b)=\{ j: \a_j(a) \ne \a_j(b) \}$.  Then
$$
\split
H_{a,b}(u) &:=\sum_{j\in J(a,b)} e^{\b_j u} \( f_j(u,\a_j(a))-
  f_j(u,\a_j(b))\) \\
&= \Delta_{a,b}(u) + o\left(e^{\beta_m u}\right)\quad(u\to\infty).
\endsplit
\tag{5.3}
$$
Lastly, define the periodic functions
$$
\split
w_{j,\a}(u) &= \Re \sum_{k=1}^2 c_{j,k} e(k\a/n_j) 
\frac{e^{ik\g u}}{\b_j+ik\g} \\
&= \sum_{k=1}^2 \frac{c_{j,k}}{\sqrt{k^2\g^2+\b_j^2}}
\sin\( k\g u+\frac{2\pi k\a}{n_j}+\tan^{-1}\frac{\b_j}{k\g} \).
\endsplit
$$
By Lemma 1.1 (the asymptotic for $f(\rho)$) and (5.2),
$$
f_j(u,\a) = w_{j,\a}(u) + O(1/u). \tag5.4
$$
Similarly,
$$
\split
\frac{d}{du} f_j(u,\a) &= w'_{j,\a}(u) + \sum_{k=1}^2 \frac{c_{j,k} 
e(\frac{k\a}{n_j})}
{\b_j+ik\g} \left[ (1-\b_j u)e^{-\b_j u}  \int_2^{e^u}
\frac{v^{\b_j+ik\g}}{v\log^2 v}\, dv + \frac{e^{ik\g u}}{u} \right] \\
&= w'_{j,\a}(u) + O(1/u).
\endsplit\tag5.5
$$

Each function $w_j(u,\a)$ is periodic in $u$ with period $2\pi/\g$.
We choose the numbers $\b_j$ and $c_{j,k}$ so that the functions
$w_{j,\a}$ have several properties:
\medskip

\item{(A)}  For each $j$ and each pair of distinct integers 
$\a_1,\a_2\in [0,n_j-1]$, the equation
$$
w_{j,\a_1}(u) = w_{j,\a_2}(u)
$$
has only two solutions in $[0,2\pi/\g)$.
Call them $\theta_v(j,\a_1,\a_2)$, $v=1,2$;

\item{(B)} All the numbers $\theta_v(j,\a_1,\a_2)$ are nonzero
and distinct, that is
$$
\theta_{v_1}(j_1,\a_1,\a_2) = \theta_{v_2}(j_2,\a_3,\a_4) \text{ implies }
v_1=v_2, j_1=j_2, \{\a_1,\a_2\} = \{\a_3,\a_4\};
$$

\item{(C)}
For all $j$, $v$ and distinct $\a_1,\a_2$, if $\theta=\theta_v(j,\a_1,\a_2)$
then $w'_{j,\a_1}(\theta)-w'_{j,\a_2}(\theta) \ne 0$;

\item{(D)}
Let $1\le j'<j\le m$, $v\in \{1,2\}$, 
distinct $\a_1,\a_2\in [0,n_{j'}-1]$.
Suppose $\a_3,\a_4,\a_5,\a_6\in [0,n_j-1]$ with
$(\a_3,\a_4) \ne (\a_5, \a_6)$ and not both $\a_3=\a_4$ and $\a_5=\a_6$.
If $\theta=\theta_v(j',\a_1,\a_2)$, then 
$$
w_{j,\a_3}(\theta)-w_{j,\a_4}(\theta) -[w_{j,\a_5}(\theta)-w_{j,\a_6}(\theta)
 ] \ne 0.
$$

\noindent
Note: some cases of (D) are redundant, being covered by property (B).
For example, if $\a_3=\a_4$ and $\a_5\ne \a_6$, or if
$\a_3=\a_5$ and $\a_4\ne \a_6$, or if $\a_3=\a_6$ and $\a_4=\a_5$.
\medskip

\noindent
For integral $\ell$ let $u_\ell=\frac{2\pi}{\g} \ell$.
We claim the following hold for large $\ell$ (depending on $\FF$).
Throughout the remainder of this proof, $o(1)$ refers to a function
of $\ell$ which tends to 0 as $\ell\to\infty$.
\medskip

\item{(i)} At $u=u_\ell$, $u_{\ell+1}, \ldots$, the ordering
 of the functions $F_{q,a}(e^u)$ $(a\in \fq)$ is the same;

\item{(ii)} For distinct $a,b\in \fq$, the sign
changes of $\Delta_{a,b}(u)$ on $[u_\ell,u_{\ell+1}]$
occur within two intervals $I_1(a,b)$ and $I_2(a,b)$.
All $\phi(q)(\phi(q)-1)$ of these intervals are disjoint,
and the sign of each function $\Delta_{c,d}(u)$ 
at the endpoints of $I_v(a,b)$ depends only on $a,b,c,d$ and $v$.
\medskip

Together, (i) and (ii) imply the theorem.  Indeed, the possible
orderings of $F_{q,a_i}(e^u)$ ($1\le i\le r$)
are precisely the orderings occurring
at the endpoints of the intervals $I_v(a_i,a_j)$.  There are 
$r(r-1)$ such intervals, and the ordering remains constant between
two such intervals, so there are at most $r(r-1)$ different orderings. 

First we prove (i).  Let $W(j,\a) = w_{j,\a}(0)$ for each $j,\a$ and
let $L \ge \ell$.  For each $a,b\in \fq$ let
$j_0=j_0(a,b):=\min \{j\in J(a,b) \}$.  By (5.3) and (5.4),
$$
\Delta_{a,b}(u_L) = \exp \( \b_{j_0} u_L \) 
\( W(j_0,\a_{j_0}(a)) - W(j_0,\a_{j_0}(b)) + o(1) \).
$$
By (B), $W(j_0,\a_{j_0}(a)) \ne W(j_0,\a_{j_0}(b))$ and so
$\Delta_{a,b}(u_L)$ has constant sign for $L\ge \ell$.

Next we prove (ii).  Throughout suppose $u_{\ell} \le u \le u_{\ell+1}$.
For sufficiently small $\delta$ (depending only on the functions $w_{j,\a}$)
let
$$
M(a,b)=\left\{ u\in [u_\ell,u_{\ell+1}]: |H_{a,b}(u)| \le
\delta e^{\beta_m u} \right\}.
$$
By (5.3), $\Delta_{a,b}(u)=0$ implies
 $u\in M(a,b)$.  Let $j_0=j_0(a,b)$, $\a_1=\a_{j_0}(a)$
and $\a_2=\a_{j_0}(b)$.
By (5.3), for $u\in M(a,b)$,
$$
|f_{j_0}(u,\a_1)-f_{j_0}(u,\a_2)| = o(1).
$$
which by (5.4) implies that for any fixed $\eta>0$, if $\ell$ is large enough,
$$
|w_{j_0,\a_1}(u) - w_{j_0,\a_2}(u)| \le \eta. \tag5.6
$$
Let $Y$ be the set of $u$ satisfying (5.6).  By (A) and (C),
if $\delta$ and $\eta$ are small enough then $Y$ is the union of
two short intervals $K_1,K_2$, where
$\theta_v({j_0},\a_1,\a_2)\in K_v$ for $v=1,2$.
By (C), for some $\e>0$,
$w'_{j_0,\a_1}(u) - w'_{j_0,\a_2}(u)$ has constant sign and is
at least $\e$ in magnitude on each interval $K_v$. 
By (5.3), (5.4) and (5.5), $H_{a,b}(u)$ is monotone
on each of $K_1$ and $K_2$.
Therefore, $I_1(a,b):=M(a,b) \cap K_1$ and $I_2(a,b):=M(a,b) \cap K_2$
are closed intervals.  At the endpoints of $I_1(a,b)$ and $I_2(a,b)$, 
$|H_{a,b}(u)| = \delta e^{\b_m u}$ and thus
$\sgn \Delta_{a,b}(u) = \sgn H_{a,b}(u)$.
For each $v$, the sign of $H'_{a,b}(u)$ on $K_v$ thus determines the sign of
$\Delta_{a,b}(u)$ at the endpoints of $I_v(a,b)$.  This in turn depends
only on
the sign of $w'_{j_0,\a_1}(u) - w'_{j_0,\a_2}(u)$ on $K_v$, which does
not depend on $\ell$.

Next, suppose $u\in I_v(a,b)$ and $\{a,b\} \ne \{c,d\}$.
Let $j_0=j_0(a,b)$, $\a_1=\a_{j_0}(a)$,
$\a_2=\a_{j_0}(b)$.  Then
$$
|u-\theta| = o(1), \quad\text{where } \theta=\theta_v(j_0,\a_1,\a_2). \tag5.7
$$
Let $j_1=j_0(c,d)$, $\a_3=\a_{j_1}(c)$, $\a_4=\a_{j_1}(d)$.  By (5.3),
(5.4), and (5.7),
$$
\Delta_{c,d}(u) = e^{\b_{j_1} u} [ w_{j_1,\a_3}(\theta) - 
w_{j_1,\a_4}(\theta) + o(1)].
$$
If $j_0\ne j_1$ or $\{\a_1,\a_2\} \ne \{\a_3,\a_4\}$,
$w_{j_1,\a_3}(\theta)- w_{j_1,\a_4}(\theta) \ne 0$ by (B),
so $\Delta_{c,d}(u)$
has constant sign depending only on $a,b,c,d,v$.
Next, suppose $j_0=j_1$ and $\{\a_1,\a_2\} = \{\a_3,\a_4\}$.  By
swapping $c$ and $d$ if necessary, we may suppose that $\a_1=\a_3$,
$\a_2=\a_4$.  Let
$$
j_2 = \min \{ j\in J(a,b) \cup J(c,d) : 
(\a_j(a),\a_j(b)) \ne (\a_j(c),\a_j(d))\}.
$$
Such $j_2$ exists because $\{a,b\} \ne \{c,d\}$.
Also, by our assumptions on $j_0,\a_1,\ldots,\a_4$, we have $j_2 > j_0$.
By (5.3), (5.4), and (5.7),
$$
\split
\Delta_{c,d}(u) &= [\Delta_{c,d}(u) - \Delta_{a,b}(u) ] + \Delta_{a,b}(u) \\
&= e^{\b_{j_2} u} \bigl( w_{j_2}(\theta,\a_{j_2}(c)) - 
  w_{j_2}(\theta,\a_{j_2}(d))
  - [  w_{j_2}(\theta,\a_{j_2}(a)) - w_{j_2}(\theta,\a_{j_2}(b)) ] \\
& \qquad + o(1) \bigr).
\endsplit
$$
By (D), the right side has constant sign, depending only on $a,b,c,d,v$.
This completes the proof of (ii).

\bigskip

It remains to select numbers $\b_1, \ldots, \b_m$ and
$c_{j,k}$ so that (A)-(D) are satisfied.  Write for short
$$
z_j=\frac{\b_j}{\g}, \quad \e_{j} = \tan^{-1} z_j, \quad \nu_j = \tan^{-1}
\frac{z_j}{2}.
$$
Let $M$ be a large integer, depending only on $q$.  We think of 
$\b_i$ and $M$ as being fixed, while $\g\to \infty$.  In what follows
constants implied by $O$ and $\ll$ will not depend on $M$ or on $\g$.

If $n_j=2$, take $c_{j,1}=1$ and $c_{j,2}=0$.  In this case, we have
$$
w_{j,\a}(u) = \frac{\sin(\g u+\pi \a + \e_j)}
{\sqrt{\g^2+\b_j^2}}. \tag5.8
$$
If $n_j \ge 4$, we take $c_{j,1}=M$, $c_{j,2}=1$.  In this case
$$
w_{j,\a}(u) = \frac{M \sin(\g u+ \frac{2\pi \a}{n_j} + \e_j)}
{\sqrt{\g^2+\b_j^2}} + \frac{\sin(2\g u+ \frac{4\pi \a}{n_j} + \nu_j)}
{\sqrt{4\g^2+\b_j^2}}.\tag5.9
$$
In particular, for fixed $j$, $w_{j,\a}(u)=w_{j,0}(u+\frac{2\pi \a}{\g n_j})$.
Let $J_1=\{j: n_j=2\}$ and $J_2=\{j:n_j\ge 4\}$.
The functions $w_{j,0}(u)$ with $j\in J_1$ are 
very close to the function $\frac{1}{\g} \sin(\g u)$,
and the functions $w_{j,0}(u)$ with
$j\in J_2$ are all very close to the function 
$\frac{M}{\g} \sin(\g u) + \frac{1}{2\g} \sin(2\g u)$.
It is important, however, that the actual functions $w_{j,0}$ $(j\in J_2)$
are not odd nor are they a shift of an odd function.

Assume throughout that $0\le u< 2\pi/\g$.  Consider first the equation 
$$
w_{j,\a_1}(u) = w_{j,\a_2}(u), \qquad \text{ where }
0 \le \a_1 < \a_2 \le n_j-1. \tag5.10
$$
If $j\in J_1$ then $\a_1=0,\a_2=1$ and the solutions of (5.10) are 
$$
\g u\in \{ \pi-\e_j,2\pi-\e_j \}. \tag5.11
$$
Since the numbers $\e_j$ are distinct and $O(1/\g)$ in magnitude,
all such solutions
(for varying $j$) are distinct and non-zero. 
Similarly, when $j\in J_2$ and $\a_2=\a_1 + \frac12 n_j$, (5.9) implies
that the solutions of (5.10) are
$$
\g u\in \{ \pi(1-2\a_1/n_j)-\e_j,  \pi(2-2\a_1/n_j)-\e_j \}. \tag5.12
$$
Again these numbers are all distinct and non-zero (for varying $j$ and
$\a_1$), and distinct from the numbers in (5.11).
Finally, suppose $j\in J_2$ and $\a_2 - \a_1 \ne \frac12 n_j$.
We make use of the following
expression for $w_{j,\a}(u)$ which avoids square roots:
$$
\g w_{j,\a}(u) = M \frac{z_j\cos(\omega)+ \sin(\omega)}{1+z_j^2} + 
\frac{z_j\cos(2\omega)+ 2 \sin(2\omega)}{4+z_j^2},\, 
\omega=\g u+\frac{2\pi\a}{n_j}.
\tag5.13
$$
Making the change of variables $y=\g u+\frac{\pi}{n_j}(\a_1+\a_2)$, define
$$
g(y) = g(y;j,\a_1,\a_2) = \frac{-\g}{2} (1+z_j^2)(4+z_j^2)
(w_{j,\a_1}(u)-w_{j,\a_2}(u)).
$$ 
Using some trigonometric identities with (5.13), we have
$$
\split
g(y) &= M (4+z_j^2) \sin B (\cos y - z_j \sin y) 
+ (1+z_j^2) \sin 2B (2\cos 2y - z_j \sin 2y) \\
&= \sin B \left[ M (4+z_j^2)(\cos y - z_j \sin y) +  2(1+z_j^2) \cos B
(2\cos 2y - z_j \sin 2y) \right],
\endsplit\tag5.14
$$
where
$$
B=\frac{\pi(\a_2-\a_1)}{n_j} \in \left\{ \frac{k\pi}{n_j} : 1\le k\le
n_j-1, k\ne n_j/2 \right\}. \tag5.15
$$
Since $\sin B \ne 0$ by (5.15), combining (5.10) and (5.14) gives
 the approximation
$$
4M\cos y + 4\cos 2y \cos B = O(M/\g). \tag5.16
$$
We may assume $\g \ge M$.
Thus $|\cos y| \ll 1/M$ and consequently $|\sin y|=1+O(1/M^2)$,
$\cos 2y=-1+O(1/M^2)$, and $|y\pm \pi/2| \ll 1/M$.  For such $y$,
$|g'(y)| \gg M$, so there are exactly two solutions of (5.10), one with $y$
 near $\pi/2$ and the other with $y$ near $-\pi/2$.
This proves (A).  When $u=0$, i.e. $y=\frac{\pi}{n_j}(\a_1+\a_2)$,
(5.15) implies $|g(y)| \gg M$ unless $\a_1+\a_2 \in \{ n_2/2, 3n_j/2\}$.
In this case $g(y)=-4\cos B\sin B + O(M/\g) \ne 0$ by (5.15).
This proves that every $\theta_v(j,\a_1,\a_2) \ne 0$.

For the second part of (B), consider the equation
$$
\theta_{v_1}(j_1,\a_1,\a_2)=\theta_{v_2}(j_2,\a_3,\a_4).
$$
This implies that for some $u$,
$$
w_{j_1,\a_1}(u)=w_{j_1,\a_2}(u), \quad w_{j_2,\a_3}(u)=w_{j_2,\a_4}(u).
\tag{5.17}
$$
We cannot have $j_1=j_2\in J_1$.  First suppose $j_1,j_2\in J_2$.
We may assume $\a_1<\a_2$, $\a_3<\a_4$ and either $j_1\ne j_2$ or
$\{\a_1,\a_2\} \ne \{\a_3,\a_4\}$.
If $j_1=j_2=j$ and $\a_i=\a_k$ for some $i\ne k$, then
$w_{j,\a_1}(u)=w_{j,\a_2}(u)=w_{j,\a_3}(u)=w_{j,\a_4}(u)$,
the set $\{\a_1,\a_2,\a_3,\a_4\}$ contains three distinct elements,
and the function $w_{j,0}$ takes some value three times,
which is impossible.

By (5.14), we have the system of equations
$$
\split
M(4+z_{j_1}^2)(\cos y_1 - z_{j_1}\sin y_1) + 2(\cos B_1)(1+z_{j_1}^2)
(2\cos 2y_1 - z_{j_1}\sin 2y_1)&=0, \\
M(4+z_{j_2}^2)(\cos y_2 - z_{j_2}\sin y_2) + 2(\cos B_2)(1+z_{j_2}^2)
(2\cos 2y_2 - z_{j_2}\sin 2y_2)&=0,
\endsplit\tag5.18
$$
where 
$$
y_1=\g u+\frac{\pi(\a_1+\a_2)}{n_{j_1}}, 
y_2 = \g u+\frac{\pi(\a_3+\a_4)}{n_{j_2}},
B_1=\frac{\pi(\a_2-\a_1)}{n_{j_1}}, B_2=\frac{\pi(\a_4-\a_3)}{n_{j_2}}.
$$
As before, $|\cos y_k|=O(1/M)$ for $k=1,2$.  Since $y_1-y_2$ is an integral
multiple of $\pi/\phi(q)$ and $\g$ is large, 
$\cos y_1 = \pm\cos y_2$.   As a consequence, $\cos 2y_1=\cos 2y_2
=-1+O(1/M^2)$ and so by (5.18),
$$
4M \cos y_1 + 4\cos 2y_1 \cos B_k = O(M/\g) \qquad (k=1,2).
$$
This in turn implies that either $\cos y_1 = \cos y_2$ and $B_1=B_2$
or that $\cos y_1 = -\cos y_2$ and $B_1=\pi-B_2$. 
Consider four cases: (i) $\cos y_1=\cos y_2$, $\sin y_1=\sin y_2$,
(ii) $\cos y_1=\cos y_2$, $\sin y_1=-\sin y_2$,\newline
(iii) $\cos y_1=-\cos y_2$, $\sin y_1=\sin y_2$,
(iv) $\cos y_1=-\cos y_2$, $\sin y_1=-\sin y_2$.  
In cases (ii) and (iii), subtracting or adding the two equations in (5.18) 
yields
$$
(z_{j_1}+z_{j_2})(4M\sin y_1 - 2\cos B_1 \sin 2y_1)=O(M/\g^2),
$$
which is not possible given that $|\sin y_1|=1+O(1/M^2)$.  In case (i)
$y_1=y_2$ and, together with $B_1=B_2$, implies that $j_1\ne j_2$
and hence $z_{j_1} \ne z_{j_2}$.  In case (iv)
$|y_1-y_2|=\pi$ and, together with $B_1=\pi-B_2$ implies that two of the 
numbers $\frac{\a_1}{n_{j_1}}$,  $\frac{\a_2}{n_{j_1}}$,  
$\frac{\a_3}{n_{j_2}}$,  $\frac{\a_4}{n_{j_2}}$ are equal, therefore,  
$j_1\ne j_2$ and $z_{j_1} \ne z_{j_2}$.  

Again subtracting the two equations in (5.18) in case (i) and adding
in case (iv) produces
$$
(z_{j_1}-z_{j_2})(4M \sin y_1 - 2\cos B_1 \sin 2y_1)= O(M/\g^2),
$$
which likewise gives a contradiction.  Therefore, (5.17) is impossible
when $j_1,j_2\in J_2$.

If $j_1\in J_1$, $j_2\in J_2$, then by (5.11), $\g u\in \{\pi-\e_{j_1}, 
2\pi-\e_{j_1}\}$. We have seen that (5.17) is impossible in the case
$\a_4=\a_3 + \frac12 n_{j_2}$. Assume that the last equality does not hold
and define $y=\g u + \pi(\a_3+\a_4)/n_{j_2}$, $B=\pi(\a_4-\a_3)/n_{j_2}$
By (5.17), $g(y;j_2,\a_3,\a_4)=0$.
From (5.14), $|\cos y|\ll 1/M$, and
thus $\pi(\a_3+\a_4)/n_{j_2}\in \{\pi/2,3\pi/2\}$.
This implies the stronger inequality $|\cos y|=|\sin \e_{j_1}| \le 1/\g$
and as a consequence $\cos 2y \le -1+2/\g^2$.  Applying (5.14)
 again we see that
$4\cos B\cos 2y = O(M/\g)$, which by (5.15) is impossible.
This completes the proof of (B).

Next we verify (C).  If $j\in J_1$, then $\a_1=0$, $\a_2=1$,
$\theta$ satisfies $\sin(\g\theta+\e_j)=0$ and thus $\cos(\g\theta+\e_j)=
\pm 1$, so $w'_{j,0}(\theta) \ne w'_{j,1}(\theta)$.  If $j\in J_2$ and
$|\a_1-\a_2|=\frac12 n_j$, the situation is the same as with $j\in J_1$
by (5.9).  Assume $0\le \a_1 < \a_2 < n_j$ and $\a_2-\a_1\ne \frac12 n_j$.
Define $B$ and $y$ as in (5.14), (5.15). 
Suppose $u$ satisfies (5.10) and also
the equation $w'_{j,\a_1}(u) = w'_{j,\a_2}(u)$.  Differentiating (5.14) gives
$$
\sin y (4M + 16\cos y \cos B)=O(M/\g),
$$
which is impossible since $|\sin y|=1+O(1/M^2)$.
Thus condition (C) is verified.

We verify condition (D) indirectly.  Condition (B) covers the situation
when $\a_3=\a_4$, $\a_5=\a_6$, $\a_3=\a_5$, $\a_4=\a_6$ or there are at
 most two distinct values among $\a_3,\ldots,\a_6$.  Henceforth assume
none of these conditions occurs and, moreover, $\a_3<\a_4$, $\a_5<\a_6$.
Fix $j',v,\a_1,\a_2$ and put $u=\theta_v(j',\a_1,\a_2)$.  
Fix $j,\a_3,\ldots,\a_6$ and define
$$
y_1=\g u+\frac{\pi(\a_3+\a_4)}{n_j}, y_2 = \g u+\frac{\pi(\a_5+\a_6)}{n_j},
B_1=\frac{\pi(\a_4-\a_3)}{n_j}, B_2=\frac{\pi(\a_6-\a_5)}{n_j}.
$$
Using (5.14), the equation in (D) becomes
$$
\multline
P(z_j):=M(4+z_j^2) \left[\sin B_1(\cos y_1- z_j\sin y_1)
- \sin B_2(\cos y_2- z_j\sin y_2)\right]\\ 
+(1+z_j^2) \left[ \sin 2B_1(2\cos 2y_1-z_j \sin 2y_1) 
- \sin 2B_2(2\cos 2y_2-z_j\sin 2y_2) \right] = 0.
\endmultline\tag5.19
$$
We shall prove that for large $M$ the polynomial $P$ is not identically zero.
The conclusion is that given $\b_1,\ldots,\b_{j-1}$, there are a finite
number of $\b_j$ which would lead to failure of (D).  Consequently, we
can always choose an admissible $\b_j$ from within a short interval.

Note that we have proved (by (5.14)) that $\g u=\frac{\pi\a}{n_{j'}}+O(1/M)$
for some integer $\a$. Therefore, there are $\tilde y_1=\frac{\pi l_1}
{\phi(q)}$ and $\tilde y_2=\frac{\pi l_2}{\phi(q)}$
with some integers $l_1$ and $l_2$ such
that $y_1=\tilde y_1+O(1/M)$, $y_2=\tilde y_2+O(1/M)$. 

Suppose that all the coefficients of $P$ are zero. The constant term 
is $4aM+O(1)$, $a=\sin B_1\cos \tilde y_1-\sin B_2\cos \tilde y_2$. 
Since $a$ can take finitely many values, we conclude from $4aM+O(1)=0$ that 
$$a=\sin B_1\cos \tilde y_1-\sin B_2\cos \tilde y_2=0.\tag5.20$$ 
In the same way, considering the coefficients of $z_j$,
we get
$$\sin B_1\sin \tilde y_1-\sin B_2\sin \tilde y_2=0.\tag5.21$$ 
Taking into account that $\sin B_1>0$, $\sin B_2>0$, we deduce from
(5.20) and (5.21) that 
$$\sin B_1=\sin B_2,\quad \tilde y_1=\tilde y_2.\tag5.22$$
Further, the last equality implies that in fact $y_1=y_2$. This in turn
implies that the sums of terms containing $M$ in the coefficients of $P$
are zero. Therefore, the conditions that the constant term and the 
coefficient of $z_j$ in $P$ are zero mean that    
$$
\cos 2y_1 (\cos B_1 - \cos B_2)=0,\quad 
\sin 2y_1 (\cos B_1 - \cos B_2)=0.
$$ 
It follows from these equalities and (5.22) that 
$$B_1=B_2.\tag5.23$$
Finally, from (5.22) and (5.23) we obtain $\a_3=\a_5$ and $\a_4=\a_6$,
which does not agree with our assumptions and completes the proof of Theorem
5.1. 
\qed\enddemo

\vfil\eject

%
%

\enddocument
\bye
\end

\enddocument
\end